\documentclass[aap]{imsart}
\RequirePackage[OT1]{fontenc}

\usepackage{chicago}

\usepackage{graphicx}
\usepackage{amssymb}
\usepackage{amsmath}
\usepackage{subfigure}

\newcommand{\qed}{\hfill $\Box$}

\newtheorem{theorem}{Theorem}
\newtheorem{thrm}{Theorem.}[section]

\newtheorem{lemma}[thrm]{Lemma}

\newcommand{\RR}{{\mathbb R}}

\newenvironment{proof}[1][Proof.]{\begin{trivlist}
\item[\hskip \labelsep {\bfseries #1}]  }{\end{trivlist}}

\begin{document}

\begin{frontmatter}

% "Title of the paper"
\title{Integral Equations and the First Passage Time of Brownian Motions}
\runtitle{Integral Equations and the FPT of Brownian Motions}

\begin{aug}
% indicate corresponding author with \corref{}
\author{\fnms{Sebastian} \snm{Jaimungal}\thanksref{t1,m1}\ead[label=e1]{sebastian.jaimungal@utoronto.ca}},
\author{\fnms{Alex} \snm{Kreinin}\thanksref{m1,m2}\ead[label=e2]{alex.kreinin@algorithmics.com}}
\and
\author{\fnms{Angelo} \snm{Valov}\thanksref{t1,m1}\ead[label=e3]{valov@utstat.utoronto.ca}}

\thankstext{t1}{This research was partly supported by grants from NSERC of Canada and the MITACS National Center of Excellence.}

\runauthor{Jaimungal, Kreinin and Valov}

\affiliation{Department of Statistics, University of Toronto\thanksmark{m1} and Algorithmics Inc.\thanksmark{m2}}

\address{Department of Statistics, University of Toronto \\ 100 St. George Street, \\  Toronto, Ontario, \\ Canada M5S 2J6 \\
\printead{e1}\\
\phantom{E-mail:\ }
\printead*{e3}}

\address{Principal Mathematician, Algorithmics Inc. \\ 185 Spadina Avenue \\ Toronto, Ontario, Canada M5T 2C6 \\  \printead{e2}}

\end{aug}

\begin{abstract}
The first passage time problem for Brownian motions hitting a barrier has been extensively studied in the literature. In particular, many incarnations of integral equations which link the density of the hitting time to the equation for the barrier itself have appeared. Most interestingly, \citeN{Peskir1} demonstrates that a master integral equation can be used to generate a countable number of new equations via differentiation or integration by parts. In this article, we generalize Peskir's results and provide a more powerful unifying framework for generating integral equations through a new class of martingales. We obtain a continuum of Volterra type integral equations of the first kind and prove uniqueness for a subclass. Furthermore, through the integral equations, we demonstrate how certain functional transforms of the boundary affect the density function. Finally, we demonstrate a fundamental connection between the Volterra integral equations and a class of Fredholm integral equations.
\end{abstract}

%\begin{keyword}[class=AMS]
%\kwd[Primary ]{}
%\kwd{}
%\kwd[; secondary ]{}

%\end{keyword}

\begin{keyword}
\kwd{First Passage Time}
\kwd{Volterra Integral Equations}
\kwd{Fredholm Integral Equations}
\kwd{Martingales}
\end{keyword}

\end{frontmatter}

\section{Introduction}
Let $(W_t)_{t\geq0}$ be a standard Brownian motion started at zero and $b:(0,\infty)\rightarrow \RR$ be a continuous function satisfying $b(0)\leq 0$. Define the first passage time (from above) of $W_t$ to the curved boundary $b(t)$ to be:
\begin{eqnarray}
\tau=\inf\{t>0; W_t \leq b(t)\} \label{def:tau}\ ,
\end{eqnarray}
with distribution function $F(t)\triangleq P(\tau \leq t)$. The first passage time (FPT) problem seeks to determine $F$ when $b$ is given, while in the inverse problem we look for $b$ given $F$. We will assume that $b(t)$ is a \textit{regular} boundary in the sense that $P(\tau=0)=0$. Sufficient conditions for regularity are given by Kolmogorov's test (see e.g. \citeN{ItoMcKean} pp. 33-35). Furthermore, we allow $b(0)=-\infty$ but we assume that whenever this is the case then there exists $\epsilon>0$ such that $b$ is monotone increasing on $(0,\epsilon]$.

The FPT problem for Brownian motion has a long history and available closed form results appear to be sparse and fragmented. The few special cases include the linear boundary, quadratic boundary (see \citeN{Salminen88} and \citeN{Groeneboom89}), and square-root boundary (see \citeN{Breiman66}, \citeN{RicciardiSacerdoteSato84} and \citeN{NovikovFrishlingKordzakhia99} among others). The celebrated method of images allows one to, at least theoretically, solve the problem for a class of boundaries, $b_a(t)$, which are solutions, for each fixed $t$, of implicit equations of the type $$\int_0^{\infty}e^{-ub_a(t)-u^2t/2}Q(du)=a\ ,$$
where $a>0$ and $Q$ is a positive $\sigma$-finite measure (see \citeN{Lerche86}). In this case the density function $f_a(t)$ is given by $$
f_a(t)=\frac{\phi(b_a(t)/\sqrt{t})}{2t^{3/2}}\frac{\int_0^{\infty}\theta\phi(\frac{b_a(t)+\theta}{\sqrt{t}})Q(d\theta)} {\int_0^{\infty}\phi(\frac{b_a(t)+\theta}{\sqrt{t}})Q(d\theta)} \ .
 $$

One way to tackle the FPT problem is to derive equations linking $b$ and $F$. This is one of the primary motivations for studying integral equations in the context of the FPT. %(see \citeN{Peskir1} and the references therein).
\citeN{Peskir1} presents a unifying approach to the integral equations of Volterra type arising from the FPT. Furthermore, the author generalizes the Volterra equations of the first kind. These equations are difficult to solve analytically but they are useful in a number of areas including the study of the small time behavior of $F$ (\citeN{Peskir2}), numerical procedures yielding approximate evaluations of $F$ (\citeN{ParkSchuurmann76}, \citeN{Smith72}, \citeN{Durbin71} among others) or closed and semi-closed form approximations (e.g. \citeN{Ferebee82} and \citeN{ParkParanjape74}). Integral equations of Fredholm type are also useful in deriving known unique integral transforms of $F$ (e.g. \citeN{Shepp67}, \citeN{Novikov81}) and expansions of the FPT density (e.g. \citeN{Daniels2000}).

The first kind Volterra or Fredholm equations mentioned above can be viewed as being of the form $E(g(W_{\tau},\tau))=g(0,0)$ and as such are a direct result of the optional sampling theorem applied to an apropriate martingale $g(W_s,s)$. Using this simple martingale result, our main aim is to present a unifying approach to the integral equations arising from the FPT  and generalize the known class of integral equations. In Section 2 we examine such classes of martingales and provide a class of integral equations which generalize all previously known Volterra integral equations of the first kind. Furthermore we examine necessary and sufficient conditions for the existence of a unique solution to a subclass of these equations. In addition, we outline a method, based on the method of images, for deriving new integral equations of Volterra type. In Section 3 we apply a similar martingale approach to derive Fredholm type equations. These equations are then used to provide concise alternative derivations of known closed form results for the linear, quadratic and square-root boundaries. Finally, we show the equivalence between the Fredholm and Volterra equations of the first kind for a particular class of boundaries.

\section{Volterra Integral Equations}

The motivation behind connecting the martingale theory and the construction of integral equations for Brownian motion is perhaps best illustrated by the following well known Volterra equation (\citeN{Peskir1}):
\begin{align}
\int_0^t \phi\left(\frac{y-b(s)}{\sqrt{t-s}}\right) \ \frac{F(ds)}{\sqrt{t-s}}= \frac{1}{\sqrt{t}}\phi\left(\frac{y}{\sqrt{t}}\right) \label{eqn:IE1}
\end{align}
where $\phi$ is the standard normal density function. The equality holds for all $y< b(t)$ for continuous regular boundaries $b$. This equation can be written as
\begin{eqnarray}
E(X_{\tau}1(\tau\leq t))=X_0 \label{OPT}
\end{eqnarray}
 where the process $X_s$ is defined as $X_s=\phi\left(\frac{y-W_s}{\sqrt{t-s}} \right)$ for fixed $t>0$. Replacing $\phi$ by $\Phi$, the standard normal cdf, produces another well known equation which holds for all $y\leq b(t)$ when $b$ is continuous. Noting that $X_s$ is a real-valued martingale for $s<t$ and that $X_t1(\tau>t))=0$ a.s., equation (\ref{OPT}) can be viewed as a product of the optional sampling theorem applied to the process $X_{s\wedge t}$ and the stopping time $\tau$.\\

 Thus, the first step is to look for a class of martingales of the form $X_s\triangleq m(W_s,s),\ s<t,$ satisfying $E(|X_{\tau}|1(\tau\leq t))=\int_0^{t}|m(b(s),s)|F(ds)<\infty$ and such that $\lim_{s\uparrow t}E(X_s1(\tau>s))=0$. Suppose such a martingale exists and take a localizing sequence of stopping times $s\wedge \tau,\ s<t$. Then, applying the optional sampling theorem to $X$ and $s\wedge \tau$ and passing to the limit $s\uparrow t$, we obtain
$$X_0=\lim_{s\uparrow t}E(X_{s\wedge \tau})=\lim_{s\uparrow t}E(X_{\tau}1(\tau \leq s))=\int_0^{t}m(b(u),u)F(du)$$
by the use of the almost sure identity $W_{\tau}=b(\tau)$. The class of functions $m$ for which the process $X_s$ satisfies the above properties is a rather large class. A subclass of positive functions $m$ can be constructed using the following classical result due to \citeN{Widder44}:

\begin{theorem}[\textbf{\citeN{Widder44}}]
\label{thm:widder}
Let $u$ be a continuous, non-negative function on $I=(0,\delta)\times \RR,\ 0<\delta\leq \infty$. The following statements are equivalent:\\
1) $u$ satisfies the diffusion equation $u_s=\frac{1}{2} \ u_{xx}$ on I and $\lim_{(s,x)\rightarrow (0,e)}u(s,x)=0$ for all $e<0$\\
2) There exists a positive $\sigma$-finite measure $Q$ on $[0,\infty)$ such that $u$ can be represented as
\begin{equation}
u(s,x)=\int_0^{\infty}\frac{1}{\sqrt{s}}\phi\left(\frac{x-\theta}{\sqrt{s}}\right)\ Q(d \theta) \ .\label{eqn:widder}
\end{equation}
\end{theorem}

Given this result, define $m(x,s)\triangleq u(t-s,y-x)$ for any $t>0$ and $y<b(t)$. Then $m$ satisfies the diffusion equation $m_s=-\frac{1}{2}m_{xx}$ using the first part of the Theorem \ref{thm:widder}. Furthermore, the process $X_s\triangleq m(W_s,s),\ s<t$ is a martingale (we can check directly, by computing the double integral, that $E(|X_s|)=X_0<\infty$ for all $s<t$). Checking the first condition, $E(|X_{\tau}|1(\tau\leq t))<\infty$, we have:
$$\int_0^{t}|m(b(s),s)|F(ds)=\int_0^{\infty}Q(d\theta)\int_0^t\frac{1}{\sqrt{t-s}}\phi\left(\frac{b(s)-(z-\theta)}{\sqrt{t-s}}\right)F(ds)=u(z,t)$$
using equation \eqref{eqn:IE1}. Furthermore, note that on the set $\{\tau>s\}$ we have $W_s>b(s)$. Take $s_0$ close enough to $t$ and such that for all $s_0<s\leq t$ we have $b(s)>y$. Such $s_0$ exists since $b$ is continuous and $b(t)>y$. Then
\begin{align*}
E(X_s1(\tau>s)) = &\int_0^{\infty}Q(d\theta)\ \left.E\left[1(\tau>s)\phi\left(\frac{y-W_s-\theta}{\sqrt{t-s}}\right)\right] \right/\sqrt{t-s}\\
\leq & \int_0^{\infty}\frac{1}{\sqrt{t-s}}\phi\left(\frac{y-b(s)-\theta}{\sqrt{t-s}}\right)\ Q(d \theta)\\
= & u(t-s,y-b(s))\ .
\end{align*}
Taking the limit $s \uparrow t$ and using the limitimg behavior of the function $u$ as given in Theorem \ref{thm:widder} above we see that $\lim_{s\uparrow t}E(X_s1(\tau>s))=0$. This gives us the Volterra equation of the first kind:
\begin{eqnarray}
u(t,y)=\int_0^tu(t-s,y-b(s))F(ds)
\end{eqnarray}
for any $y<b(t),\ t>0$.

The integral representation of the function $u$ (equation \eqref{eqn:widder}) is computable for several specific ``degenerate'' cases, such as when $Q(d\theta)$ is a sum of Dirac measures or a uniform measure over a compact domain. However, we have found one other general class of measures which lead to tractable forms for $u$ itself, specifically when $Q(d\theta)=\theta ^{-p-1}d\theta$ for $p<0$. In this case by direct calculation (see \eqref{pcfintegral}) we have
$$u(t,x;p)=e^{-x^2/(4t)}D_p(-x/\sqrt{t})\sqrt{t}^{-p-1}\gamma(p)/\sqrt{2\pi}$$
where $D_p$ is the parabolic cylinder function (see Section \ref{pcf}). Note that even for $p\geq0$, this particular $u(t,x;p)$ still satisfies the diffusion equation $u_t=\frac{1}{2}u_{xx}$; furthermore,  $u(t,x;0)\sim \phi(x/\sqrt{t})/\sqrt{t}$ and $u(t,x;-1)\sim\Phi(x/\sqrt{t})$ which are the kernels of the two well known Voltera equations. These observations motivate us to examine the function
\begin{eqnarray}
m(s,x;p)=\frac{e^{-\frac{(x-y)^2}{4(t-s)}}D_p((x-y)/\sqrt{t-s})}{(t-s)^{(p+1)/2}} , \quad p,y,x \in \RR,\  \label{fn:cylinder}
\end{eqnarray}
(for a fixed $t>0$) more closely\footnote{We first came across this function through an alternative route prior to realizing the connection to the Widder's (1944) result. In fact, it is not apparent how Widder's result applies when $p>0$. However, $D_p \ (p>0)$ can be written as a linear combination of parabolic cylinder functions with $p<0$, the coefficients of which are space and time dependent.}. We now proceed to show that $m(s,W_s),\ s<t,$ is an honest martingale and derive a Volterra equation with kernel $m(s,b(s);p)$ by applying the optional sampling theorem.

Define the process $X_{s}\triangleq m(s,W_{s};p)$, the stopping time $\tau_t\triangleq\tau \wedge t$, and fix $y \in (-\infty,b(t))$. We will use the optional sampling theorem on $X_s,\ s<t,$ and $\tau_t$ -- it is important to point out that here time flows with $s$, while $t$ represents a fixed time point.

First we show that $\{X_s\}_{s<t}$ is a martingale. Using the second order differential equation \eqref{diffeqn1}, to which $D_p$ is a solution, it is straightforward to show that
\begin{eqnarray}
m_s=-\frac{1}{2}m_{xx} \ . %=\frac{e^{-z^2/4}}{(t-s)^{(p+3)/2}}\{D_{p}(z)(v+1-z^2/4)+D'_{p}(z)z/2\}
\end{eqnarray}
To check the integrability conditions, consider $E\left[|m(s,W_{s};p)|1(|W_s|>a)\right]$, $a\gg y$, $s<t$. Using the asymptotic behavior of the parabolic cylinder function (see \eqref{pcfassymptotic+} and \eqref{pcfassymptotic-}) we obtain:
\begin{align*}
& E\left[|m(s,W_{s};p)|1(|W_s|>a)\right]\\
& \hspace{1cm} = \int_{-\infty}^{-a}|m(s,x;p)|\frac{e^{-x^2/(2s)}}{\sqrt{2\pi s}}dx+\int_{a}^{\infty}|m(s,x;p)|\frac{e^{-x^2/(2s)}}{\sqrt{2\pi s}}dx\\
& \hspace{1cm} \sim \int_{-\infty}^{-a}\frac{e^{-\frac{x^2}{2s}}|(x-y)^{-p-1}|}{\sqrt{2\pi s}}dx +\int_{a}^{\infty}\frac{e^{-\frac{(x-y)^2}{2(t-s)}-x^2/(2s)}(x-y)^p}{(t-s)^{p+1/2}\sqrt{2\pi s}}dx <\infty
\end{align*}
Furthermore, $m(s,x;p)$ is a continuous function in $x$ on $[-a,a]$. Thus,
\begin{align*}
E(|X_s|) =& E|m(s,W_{s};p)| \\
=& E\left[|m(s,W_{s};p)|1(|W_s|>a)\right]+E\left[|m(s,W_{s};p)|1(|W_s|\leq a)\right]<\infty.
\end{align*}
 Therefore, $X_s$ is a martingale for all $p,y\in \RR$. For $s<t$ the process is a real valued martingale while for $s>t$ it is a complex valued martingale.

As before, on the set $\{\tau_t>s\}$, we have $W_s>b(s)$ which implies $(W_s-y)/\sqrt{t-s}>(b(s)-y)/\sqrt{t-s}\rightarrow \infty$ as $s\uparrow t$ because of the continuity of $b(.)$ and the condition $y<b(t)$. Thus, choosing $s_0$ close enough to $t$ and such that for all $s_0<s\leq t$ we have $b(s)>y$ and using the asymptotic behavior of the parabolic cylinder function (\ref{pcfassymptotic+}), we obtain
\begin{eqnarray*}
|X_s|1(\tau_t>s)&=&m(s,W_s;p)1(\tau_t>s)\\
&=&1(\tau_t>s)\frac{e^{-\frac{(W_s-y)^2}{4(t-s)}}|D_p((W_s-y)/\sqrt{t-s})|}{(t-s)^{(p+1)/2}}\\
&\sim& 1(\tau_t>s)\frac{e^{-\frac{(W_s-y)^2}{2(t-s)}}(|W_s-y|/\sqrt{t-s})^{p}}{(t-s)^{(p+1)/2}}\\
&\leq& \frac{e^{-\frac{(b(s)-y)^2}{2(t-s)}}}{(t-s)^{(2p+1)/2}}1(\tau_t>s)|W_s-y|^p
\end{eqnarray*}
In particular, $m(t,b(t);p)=\lim_{s\uparrow t}m(s,b(t);p)=0$, for all $p$ since $b(t)>y$. Furthermore,
\begin{eqnarray*}
E(1(\tau_t>s)|W_s-y|^p )&\leq& \int_{b(s)-y}^{\infty}x^p\frac{e^{-(x-y)^2/(2s)}}{\sqrt{2\pi s}}dx\\
&\leq& \frac{1}{\sqrt{2\pi s}}\int_{b(s)-y}^{\infty}x^pe^{-(x-y)^2/(2t)}dx < \infty \ \forall s\leq t\ .
\end{eqnarray*}
Therefore,
\begin{eqnarray*}
\lim_{s\uparrow t}E|X_s|1(\tau_t>s)\leq C \lim_{s\uparrow t}\frac{e^{-\frac{(b(s)-y)^2}{2(t-s)}}}{(t-s)^{(2p+1)/2}}=0
\end{eqnarray*}
since $(b(s)-y)\rightarrow (b(t)-y)>0$. Furthermore, whenever $\tau>t$, then $|X_{\tau_t}|=|X_t|=\lim_{s\uparrow t}|X_s|=0$ since $\tau> t$ implies $(W_s-y)>(b(s)-y)$ for all $s<t$ and $b(t)-y>0$. Therefore,
\begin{align*}
E(|X_{\tau_t}|) &= E(|X_{\tau}|1(0<\tau\leq t))\\
&= E(|m(\tau,b(\tau);p)|1(\tau\leq t))=\int_0^t|m(s,b(s);p)|F(ds) \ .
\end{align*}
We already have that $m(t,b(t);p)=0$ for all $p$. Since $m$ is a continuous function (since $b$ is continuous) it follows that the last integral above is finite provided that $$\int_0^{\epsilon}|m(s,b(s);p)|F(ds)<\infty$$ for some small positive $\epsilon$. This is the case when $b(0)>-\infty$ since $m(0,b(0);p)<\infty$, so let us assume that $b(0)=-\infty$. Choosing $\epsilon$ small enough and using the asymptotic behavior of the parabolic cylinder function we have
\begin{eqnarray}
\int_0^{\epsilon}|m(s,b(s);p)|F(ds)\sim \int_0^{\epsilon}|(b(s)-y)^{-p-1}|F(ds) \label{eqn:b(0)}
\end{eqnarray}
and for $p> -1$ the last integral is finite since $(b(s)-y)^{-p-1}\rightarrow 0$ as $s\downarrow 0$ and equals $F(\epsilon)$ for $p=-1$. The case $p<-1$ follows from Lemma \ref{appendixlemma} and (\ref{eqn:finite}) in Appendix A. Therefore, for all $p\in \mathbb \RR$ and $y<b(t)$, we have $E|X_{\tau_t}|<\infty$ and by the optional sampling theorem $X_0=E(X_{\tau_t})$; we have then proved the following result.
\begin{theorem}
Let $(W_t)_{t\geq0}$ be a standard Brownian motion with $W_0=0$. Let $b:(0,\infty)\mapsto \RR$ be a continuous function satisfying $b(0)\leq 0$. Let $\tau$ be the first-passage time of $W$ to $b$, and let $F$ denote its distribution function. Then for all $p\in \RR$ and $y<b(t)$ the following system of integral equations is satisfied:
\begin{eqnarray}
\frac{e^{-\frac{y^2}{4t}}D_{p}(-y/\sqrt{t})}{t^{(p+1)/2}}=\int_{0}^te^{-\frac{(b(s)-y)^2}{4(t-s)}}\frac{D_{p}((b(s)-y)/\sqrt{t-s})}{(t-s)^{(p+1)/2}}F(ds) \label{eqn:general}
\end{eqnarray}
where $F$ is the distribution of $\tau$.
\end{theorem}

The set of integral equations (\eqref{eqn:general}) reduce to a set of well known integral equations when $p=n$, a non-negative integer, in which case \eqref{eqn:general} becomes
\begin{align}
\frac{e^{-\frac{y^2}{2t}}H_n(-y/\sqrt{2t})}{t^{(n+1)/2}}=\int_{0}^te^{-\frac{(b(s)-y)^2}{2(t-s)}}\frac{H_n((b(s)-y)/\sqrt{2(t-s)})}{(t-s)^{(n+1)/2}}F(ds)  \ . \label{eqn:hermite}
\end{align}
Here, $H_n$ are the Hermite polynomials of degree $n$ (see (\ref{pcfHn})). These equations were derived in e.g. \citeN{Peskir1} among others. In the next section we examine the limit $y\uparrow b(t)$ which allows the density and boundary to be tightly bound via the integral equations without the appearance of the arbitrary parameter $y$. Afterwards, we provide a richer class of examples.

\subsection{Passage to the limit}

The next step is to investigate what conditions on the boundary $b$ are necessary to allow the limit $y\uparrow b(t)$ in (\ref{eqn:general}) to be taken. This limit is not straightforward for all values of the parameter $p$. To see this let us compute the limit as $y\uparrow b(t)$ in equation (\ref{eqn:general}) with $p=n=1$ assuming $b(t)$ is continuously differentiable on $(0,\infty)$ and $b(0)<0$.

First, in this case, for $t_0>0$ there exists some $\epsilon>0$ such that $\epsilon \leq e^{-\frac{(b(s)-b(t))^2}{2(t-s)}}$ for all $t_0\leq s \leq t$ since $\lim_{s\uparrow t}\frac{b(s)-b(t)}{\sqrt{t-s}}=b'(t).0=0$. Then we have
\begin{eqnarray*}
\epsilon \int_{t_0}^t\frac{F(ds)}{\sqrt{t-s}}&\leq&
\int_{t_0}^t\frac{e^{-\frac{(b(s)-b(t))^2}{2(t-s)}}}{\sqrt{t-s}}F(ds)
< \int_{0}^t\frac{e^{-\frac{(b(s)-b(t))^2}{2(t-s)}}}{\sqrt{t-s}}F(ds)\\
&=& \int_{0}^t \liminf_{y\uparrow b(t)}\frac{e^{-\frac{(b(s)-y)^2}{2(t-s)}}}{\sqrt{t-s}}F(ds)\leq \liminf_{y\uparrow b(t)}\int_{0}^t\frac{e^{-\frac{(b(s)-y)^2}{2(t-s)}}}{\sqrt{t-s}}F(ds)\\
&=&\frac{e^{-\frac{b^2(t)}{2t}}}{\sqrt{t}}<\infty
\end{eqnarray*}
where the last equality follows from (\ref{eqn:hermite}) with $n=0$. Thus, when $b(t)$ is differentiable $\int_{0}^t\frac{F(ds)}{\sqrt{t-s}}<\infty$ and therefore $\int_0^t\frac{|b(t)-b(s)|}{(t-s)^{3/2}}F(ds)<\infty$ since $\frac{|b(t)-b(s)|}{t-s}$ is finite for all $0<s\leq t$ and in the neighborhood of $0$ the finiteness follows from Lemma \ref{appendixlemma} and (\ref{eqn:finite}).

Second, for such boundaries, the corresponding density function of $\tau$ is continuous i.e. $F(ds)=f(s)ds$ where $f$ is continuous on $[0,\infty)$ and $f(0)=0$ (see \cite{Peskir1} and \cite{Peskir2}). As a result,
\begin{eqnarray*}
\frac{e^{-\frac{b^2(t)}{2t}}b(t)}{t^{3/2}}&=&\lim_{y\uparrow b(t)}\int_{0}^te^{-\frac{(b(s)-y)^2}{2(t-s)}}\frac{(b(s)-y)}{(t-s)^{3/2}}F(ds)\\
&=&\lim_{y\uparrow b(t)}\int_{0}^te^{-\frac{(b(s)-y)^2}{2(t-s)}}\frac{(b(s)-b(t))}{(t-s)^{3/2}}F(ds)+\lim_{y\uparrow b(t)}\int_{0}^te^{-\frac{(b(s)-y)^2}{2(t-s)}}\frac{(b(t)-y)}{(t-s)^{3/2}}F(ds)\\
&=&\int_{0}^te^{-\frac{(b(s)-b(t))^2}{2(t-s)}}\frac{(b(s)-b(t))}{(t-s)^{3/2}}F(ds)\\
&+&\lim_{z\downarrow 0}2\int_{0}^{\infty}1(u\geq z)\exp\left(-\frac{1}{2}u^2\left(1+\frac{b(t-tz^2/u^2)-b(t)}{z\sqrt{t}}\right)^2\right)f(t-tz^2/u^2)du
\end{eqnarray*}
where we have used the substitutions $u=\frac{b(t)-y}{\sqrt{t-s}}$ and $z=\frac{b(t)-y}{\sqrt{t}}$ in the third equality above. For large $u\gg z$, $\frac{b(t-tz^2/u^2)-b(t)}{z\sqrt{t}}\approx 0$ and thus there exists a positive constant $a<1$ such that $\exp\left(-\frac{1}{2} u^2(1+\frac{b(t-tz^2/u^2)-b(t)}{z\sqrt{t}})^2\right)\leq e^{-au^2/2}$ for $u\gg z$. Therefore, since $f$ is uniformly bounded, by the dominated convergence theorem we obtain
\begin{eqnarray*}
\frac{e^{-\frac{b^2(t)}{2t}}b(t)}{t^{3/2}}&=&\int_{0}^te^{-\frac{(b(s)-b(t))^2}{2(t-s)}}\frac{(b(s)-b(t))}{(t-s)^{3/2}}f(s)ds\\
&+&2\int_{0}^{\infty}\lim_{z\downarrow 0}1(u\geq z)\exp\left(-\frac{1}{2}u^2\left(1+\frac{b(t-tz^2/u^2)-b(t)}{z\sqrt{t}}\right)^2\right)f(t-tz^2/u^2)du\\
&=&\int_{0}^te^{-\frac{(b(s)-b(t))^2}{2(t-s)}}\frac{(b(s)-b(t))}{(t-s)^{3/2}}f(s)ds+\sqrt{2\pi}f(t)
\end{eqnarray*}
since $\lim_{z\downarrow 0}\frac{b(t-tz^2/u^2)-b(t)}{z\sqrt{t}}=0$. This last equality can be rewritten as
\begin{eqnarray}
\frac{\phi(b(t)/\sqrt{t})b(t)}{t^{3/2}}=f(t)+\int_{0}^t\phi \left(\frac{b(t)-b(s)}{\sqrt{t-s}}\right)\frac{(b(s)-b(t))}{(t-s)^{3/2}}f(s)ds\label{2kindvolterra} .
\end{eqnarray}

The above equation was first derived by \citeN{RicciardiSacerdoteSato84} (see also \citeN{Peskir1}). It demonstrates the complexity involved in exchanging the limit (as $y\uparrow b(t)$) and the integral in our new class of integral equations (\ref{eqn:general}) -- even for the ``simple'' case of $p=n=1$. Nonetheless, we are able to compute this limiting case for a subclass of integral equations and the next result provides the required conditions on the boundary.
\newtheorem{corrolary}{Corrolary}
\begin{corrolary}
\label{c1}
Let $(W_t)_{t\geq0}$ be a standard Brownian motion with $W_0=0$. Let $b:(0,\infty)\mapsto \RR$ be a regular boundary and let $\tau$ be the first-passage time of $W$ below $b$, and let $F$ denote its distribution function. Then, for all $t>0$, the following system of integral equations is satisfied:
\begin{eqnarray}
\frac{e^{-\frac{b(t)^2}{4t}}D_{p}(-b(t)/\sqrt{t})}{t^{(p+1)/2}}=\int_{0}^te^{-\frac{(b(s)-b(t))^2}{4(t-s)}}\frac{D_{p}((b(s)-b(t))/\sqrt{t-s})}{(t-s)^{(p+1)/2}}F(ds) \label{eqn:generalb(t)}
\end{eqnarray}
i) For all $p\leq -1$ when $b$ is continuous on $(0,\infty)$\\
ii)For all $-1<p\leq0$ when $b$ is differentiable on $(0,\infty)$\\
iii)For all $0<p<1$ when $b$ is continuously differentiable on $(0,\infty)$
\end{corrolary}
\begin{proof} Note that $D_p(x)>0$ for all $p\leq 0$. Define $k(t)=\lim_{s\uparrow t}\frac{b(s)-b(t)}{\sqrt{t-s}}$ and $g(s;t,y)=e^{-\frac{(b(s)-y)^2}{4(t-s)}}D_{p}((b(s)-y)/\sqrt{t-s})$. The function $g$ is a continuous function in $s$ on $0<s<t$ for all $t>0$ and $y\leq b(t)$. Thus in order to apply the dominated convergence theorem we will show that $g$ is dominated by an integrable function near $s=0$ and that $g$ is finite at $s=t$ for all $y\leq b(t)$. First note that when $b(0)$ is finite then $|g(0;t,y)|$ exists for all $p$ and $y\leq b(t)$ and when $b(0)=-\infty$ then $\int_0^{\epsilon(p)}\frac{|g(s;t,y)|}{(t-s)^{(p+1)/2}}F(ds)\sim \int_0^{\epsilon(p)}(b(s)-y)^{-p-1}F(ds)<\infty$ for some $\epsilon(p)>0$ and all $p,\ y\leq b(t)$. The finiteness of the last integral follows from the fact that the integrand $(b(s)-y)^{-p-1}$ is a monotone continuous function in $y$ and thus for some $y_*$ near $b(t)$ it is dominated by $(b(s)-y_*)^{-p-1}$ which is integrable on $(0,\epsilon(p)]$ by Lemma \ref{appendixlemma} . Thus we only need to show $\lim_{s\uparrow t}\frac{g(s;t,b(t))}{(t-s)^{(p+1)/2}}<\infty$ in order to apply the dominated convergence theorem since $\lim_{s\uparrow t}\frac{g(s;t,y)}{(t-s)^{(p+1)/2}}=0$ for $y<b(t)$.\\

i) Since $\lim_{s\uparrow t}(t-s)^{-(p+1)/2}=0$ the case $|k(t)|<\infty$ is straightforward. Suppose $k(t)=\infty$. Then for $s$ close to $t$, $g(s;t,b(t))\sim e^{-\frac{(b(s)-b(t))^2}{2(t-s)}}(\frac{b(s)-b(t)}{\sqrt{t-s}})^p \rightarrow 0$ using the asymptotic behavior of $D_p(x)$ for large $x$ (see \eqref{pcfassymptotic+}). Similarly, suppose $k(t)=-\infty$ then the asymptotic behavior of $g(s;t,b(t))(t-s)^{-(p+1)/2}$ is $g(s;t,b(t))(t-s)^{-(p+1)/2}\sim (-\frac{b(s)-b(t)}{\sqrt{t-s}})^{-1-p}(t-s)^{-(p+1)/2}=(b(t)-b(s))^{-1-p}\downarrow 0$ as $s\uparrow t$ since $-p-1\geq0$ and $b$ is continuous. Therefore, taking the limit $y\uparrow b(t)$ in (\ref{eqn:general}), by the dominated convergence theorem the result follows.\\

ii) We showed that when $b$ is differentiable (and thus continuous) then $\int_{t_0}^{t}\frac{F(ds)}{\sqrt{t-s}}<\infty,\ t_0>0$. Furthermore, differentiability implies $k(t)=0$. Similarly as in part i) we see that $\lim_{s\uparrow t}g(s;t,y)/(t-s)^{p/2}=0$ for all $y\leq b(t)$ and and thus $g(s;t,y)/(t-s)^{p/2}$ is bounded on $[t_0,t]$. By the dominated convergence theorem we can exhange the limit and the integral in \ref{eqn:generalb(t)}. \\

iii) When $b$ is  continuously differentiable on $(0,\infty)$ then $f$ is continuous on $(0,t]$ for all $t>0$ \cite{Peskir1} and so $\int_0^t\frac{f(s)}{(t-s)^{(p+1)/2}}ds<\infty$ since $0<(p+1)/2<1$. Furthermore, $|g(s;t,y)|$ is bounded for $0<s\leq t$ for all $y\leq b(t)$ since $k(t)=0$. The result follows by the dominated convergence theorem. $\boxdot$
\end{proof}

Note that the differentiability condition on the boundary in part ii) can be relaxed to $|k(t)|<\infty$ for all $t>0$. In this case we still have $\int_{t_0}^t\frac{F(ds)}{\sqrt{t-s}}<\infty,\ t_0>0,$ using the same argument as before and the proof of part ii) is still valid. Also, it would be straightforward to extend the class of equations (\ref{eqn:general}) and (\ref{eqn:generalb(t)}) to the class of equations with a complex valued parameter $p$.\\

%{\bf THE NEXT PARAGRAPH NEEDS TO BE PARAPHRASED BEFORE CARRYING OUT THE ANALYSIS WHAT IS THE POINT: obtain some equations from others via a "Master equation"?}
%Let us denote the class of equations (\ref{eqn:general}) by $\{A_p(y,t)\}_{p\in \RR},\ y<b(t),$ and the class (\ref{eqn:generalb(t)}) by $\{B_p(t)\}_{p<1}$. Write $y=z-\theta,\ z<b(t),\ \theta>0,$ and let $Q(\theta)$ be a $\sigma$-finite positive measure on $[0,\infty)$. We saw above that when $Q(d\theta)=\theta^{-p-1}d\theta,\ p<0$ we can obtain equations $\{A_p(z,t)\}_{p<0}$ from $A_0(z-\theta,t)$ by applying the integral transform defined in Theorem \ref{thm:widder} w.r.t. the measure $Q$. Furthermore, in the same way we can obtain equations $\{B_p(t)\}_{p\leq-1}$ from equation $A_0(b(t)-\theta,t)$ since $Q(0)<\infty$ for $p\leq -1$. In both cases it is sufficient to assume that the boundary $b$ is continuous for equation $A_0(z-\theta,t),\ z\leq b(t),\ \theta>0$ to hold for all $t>0$. Thus, the first part of the above corollary can be obtain by simple integration without passing to the limit $y\uparrow b(t)$. Moreover, by the same integration terchnique, we can obtain equation $\{B_p(t)\}_{-1<p<0}$ from $B_0(t)$ under the hypothesis of Corollary \ref{c1}. \\

\subsection{Special Cases}
For different values of $p$ the parabolic cylinder function, $D_p$, can be written in terms of other special functions. The case when $p$ is a negative integer covers the system of equations derived in \cite{Peskir1} as we will see in \textbf{Case 3} below. Furthermore, equations (\ref{eqn:general}) and (\ref{eqn:generalb(t)}) can written in terms of the Whittaker function (see (\ref{pcfW})) or confluent hypergeometric functions using the representation of the parabolic cylinder function for all values of $p$. When $p$ is non-negative integer we already saw the connection with the Hermite polynomials which can be written in terms of the Laguerre polynomials. For $p=-1/2$ there is also a connection with the modified Bessel function of the third kind $K_{\nu}$ (see \textbf{Case 4} below). In this section we explore some of these special cases. \\

\textbf{Case 1}: $p =0$

In this case (\ref{eqn:generalb(t)}) becomes
\begin{eqnarray*}
\int_0^t \frac{e^{\frac{-(b(t)-b(s))^2}{2(t-s)}}}{\sqrt{t-s}}F(ds)=\frac{e^{-b(t)^2/2t}}{\sqrt{t}}
\end{eqnarray*}
which can be written as
\begin{eqnarray}
\int_0^t\frac{1}{\sqrt{t-s}} \phi\left(\frac{b(t)-b(s)}{\sqrt{t-s}}\right)F(ds)=\frac{1}{\sqrt{t}}\phi(b(t)/\sqrt{t}) \label{eqn:case1}
\end{eqnarray}
This equation was derived in \cite{Durbin71} who uses a previous result by Fortet(1943). \cite{Durbin71} uses the equation to obtain a numerical solution by approximating the boundary by straight line segments on subintervals $(s,s+ds)$ and using available results for crossing probabilities for linear boundaries. Subsequently  \cite{Smith72} recognizes (\ref{eqn:case1}) as a Generalized Abel equation and proposes Abel's linear transformation $T:g\rightarrow \int_0^yg(t)/\sqrt{y-t}dt$, to deal with the singularity of the kernel at $(s=t)$. He then solves the equation numerically using standard quadrature methods.\\

\textbf{Case 2}: $p=-1$

In this case (\ref{eqn:generalb(t)}) becomes
\begin{eqnarray}
\int_0^t \Phi\left(\frac{b(t)-b(s)}{\sqrt{t-s}}\right)F(ds)=\Phi(b(t)/\sqrt{t}) \label{eqn:case2}
\end{eqnarray}
Equation (\ref{eqn:case2}) was used in \cite{ParkSchuurmann76} as a basis for numerical computation of the unknown density $f$ using the idea of Volterra to discretize the equation and solve the resulting system. This equation is especially attractive for numerical computations of $f$ when $b$ is given since the kernel $K(t,s)\triangleq\Phi((b(t)-b(s))/\sqrt{t-s})$ is nonsingular in the sence that it is bounded for all $0\leq	s\leq t$. When $b(t)=c<0$ then (\ref{eqn:case2}) reads $P(\tau\leq t)=2\Phi(c/\sqrt{t})$ which is the \textit{reflection principle} for Brownian motion.\\

\textbf{Case 3}: $p=-n,\ n=1,2,3...$
In this case (\ref{eqn:generalb(t)}) becomes
\begin{eqnarray}
\int_0^t\frac{ e^{\frac{-(b(t)-b(s))^2}{4(t-s)}}}{\sqrt{2 \pi}}D_{-n}\left(\frac{(b(s)-b(t))}{\sqrt{(t-s)}}\right)(t-s)^{\frac{n-1}{2}}F(ds)=\frac{e^{-b(t)^2/4t}}{\sqrt{2\pi}}D_{-n}\left(-\frac{b(t)}{\sqrt{t}}\right)t^{\frac{n-1}{2}} \label{eqn:case3}
\end{eqnarray}

We claim that (\ref{eqn:case3}) is equivalent to Peskir's system of equations \cite{Peskir1}. Consider the kernel of the integral equation (\ref{eqn:case3}), which is of the form $\frac{1}{\sqrt{2 \pi}}e^{-x^2/4}D_{-n}(-x)=:G_{n}(x)$. Using (\ref{eqn:pcfproperty3}) we have
$$ \frac{d}{dx}G_{n+1}(x)=G_n(x) $$
and thus $$G_{n+1}(x)=\int_{-\infty}^xG_n(u)du +C$$ Taking $x=0$ and using (\ref{eqn:pcfproperty2}) we see that $C=0$. Therefore we can rewrite (\ref{eqn:case3}) as
\begin{eqnarray}
\int_0^t G_n\left(\frac{b(t)-b(s)}{\sqrt{(t-s)}}\right)(t-s)^{(n-1)/2}F(ds)=G_n(b(t)/\sqrt{t})t^{(n-1)/2} \label{eqn:peskir}
\end{eqnarray}
where $n=0,1,2,...$ and $G_{n}$ satisfies the recursion formula $G_{n+1}(x)=\int_{-\infty}^xG_n(u)du$ with $G_0(x)=\phi(x)$ since $D_0(-x)=e^{-x^2/4}$. Therefore, the system of integral equations (\ref{eqn:case3}) is equivalent to the system of equations (\ref{eqn:peskir}) which was derived in \citeN{Peskir1}. This completes the proof of the above claim.

The next two cases provide two new integral equations arising as specific cases of our general class:

\textbf{Case 4}: $p=-1/2$. In this case, using (\ref{pcfK}), (\ref{eqn:generalb(t)}) becomes
\begin{align}
& \sqrt{\frac{b(t)}{t}}e^{-\frac{b^2(t)}{4t}}K_{1/4}\left(\frac{b^2(t)}{4t}\right)\nonumber\\
&\hspace{1cm} =\int_0^t\sqrt{\frac{b(t)-b(s)}{t-s}}e^{-\frac{(b(t)-b(s))^2}{4(t-s)}}K_{1/4}\left(\frac{(b(t)-b(s))^2}{4(t-s)}\right)F(ds) \ . \label{eqn:case4}
\end{align}

\textbf{Case 5}: A new class of equations can be derived from (\ref{eqn:generalb(t)}) using the recursive relation property (\ref{eqn:pcfproperty1}) of the parabolic cylinder function. Using this relation and the class (\ref{eqn:generalb(t)}) for $p\leq -1$, we obtain:
\begin{eqnarray*}
&&\frac{e^{-\frac{b(t)^2}{4t}}}{t^{p/2}}\left\{D_{p+1}(-b(t)/\sqrt{t})+\frac{b(t)}{\sqrt{t}}D_p(-b(t)/\sqrt{t})\right\}\\
&=&\int_{0}^t\frac{e^{-\frac{(b(s)-b(t))^2}{4(t-s)}}}{(t-s)^{p/2}}\left\{D_{p+1}\left(\frac{b(s)-b(t)}{\sqrt{t-s}}\right)-\frac{b(s)-b(t)}{\sqrt{t-s}}D_p\left(\frac{b(s)-b(t)}{\sqrt{t-s}}\right)\right\}F(ds)\\
&=&\int_{0}^t(t-s)\frac{e^{-\frac{(b(s)-b(t))^2}{4(t-s)}}}{(t-s)^{(p+2)/2}}D_{p+1}\left(\frac{b(s)-b(t)}{\sqrt{t-s}}\right)F(ds)\\
&-&\int_{0}^t\frac{e^{-\frac{(b(s)-b(t))^2}{4(t-s)}}}{(t-s)^{(p+1)/2}}(b(s)-b(t))D_p\left(\frac{b(s)-b(t)}{\sqrt{t-s}}\right)F(ds)\\
&=&\frac{e^{-\frac{b(t)^2}{4t}}}{t^{p/2}}\left\{D_{p+1}(-b(t)/\sqrt{t})+\frac{b(t)}{\sqrt{t}}D_p(-b(t)/\sqrt{t})\right\}\\
&-&\int_{0}^t\frac{e^{-\frac{(b(s)-b(t))^2}{4(t-s)}}}{(t-s)^{(p+1)/2}}\left\{\frac{s}{\sqrt{t-s}}D_{p+1}\left(\frac{b(s)-b(t)}{\sqrt{t-s}}\right)+b(s)D_p\left(\frac{b(s)-b(t)}{\sqrt{t-s}}\right)\right\}F(ds)
\end{eqnarray*}
Thus, from the last equality, we derive the class of equations for $p\leq -1$:
\begin{eqnarray}
\int_{0}^t\frac{e^{-\frac{(b(s)-b(t))^2}{4(t-s)}}}{(t-s)^{(p+1)/2}}\left\{\frac{s}{\sqrt{t-s}}D_{p+1}\left(\frac{b(s)-b(t)}{\sqrt{t-s}}\right)+b(s)D_p\left(\frac{b(s)-b(t)}{\sqrt{t-s}}\right)\right\}F(ds)=0 \label{eqn:recursive}
\end{eqnarray}
For example, in the case $p=-2$ and using the results for $p=0$ and $p=-1$ together with (\ref{Dm2}), (\ref{eqn:recursive}) becomes
\begin{eqnarray}
\label{eqn:case5}
\int_0^t\left[\frac{s}{\sqrt{t-s}}\phi \left(\frac{b(t)-b(s)}{\sqrt{t-s}}\right)+b(s)\Phi \left(\frac{b(t)-b(s)}{\sqrt{t-s}} \right) \right]F(ds)=0
\end{eqnarray}\

\subsection{Uniqueness of a solution}
Next we examine sufficient conditions for the boundary $b$ such that the class of integral equations (\ref{eqn:generalb(t)}), which we denote as $\{B_p\}_{p<1}$, has a unique continuous solution. We will first investigate the uniqueness of this system of equations for the case $-1<p<1$ and then generalize to the case $p\leq-1$. Suppose that $b$ is continuously differentiable on $(0,T]$ and assume $\lim_{t\downarrow 0}|b'(t)|t^{\epsilon}<\infty$ for some $\epsilon<1/2$. Note that $\lim_{t\downarrow 0}|b'(t)|t^{\epsilon}<\infty$ implies $-\infty<b(0)<0$ (since $b$ is a regular boundary and $\epsilon<1/2$) and therefore the hitting density $f(0)=0$ \cite{Peskir2}. Therefore $F(ds)=f(s)ds$ where $f$ is continuous on $[0,\infty)$, . Denote $(p+1)/2=\lambda$ so that $0<\lambda<1$. Let
\begin{eqnarray*}
&&g_{2\lambda+1}(t)=e^{-\frac{b(t)^2}{4t}}D_{2\lambda+1}(-b(t)/\sqrt{t})/t^{\lambda}\\ &&K_{2\lambda+1}(t,s)=e^{-\frac{(b(s)-b(t))^2}{4(t-s)}}D_{2\lambda+1}((b(s)-b(t))/\sqrt{t-s})
\end{eqnarray*}
Then the class $\{B_{2\lambda+1}\}_{0<\lambda<1}$ of integral equations becomes:
\begin{eqnarray}
g_{2\lambda+1}(t)=\int_0^t\frac{K_{2\lambda+1}(t,s)}{(t-s)^{\lambda}}f(s)ds \label{eqn:unique}
\end{eqnarray}
Equations of this kind are also known as generalized Abel equations of the first kind. We know that the above equation has a continuous solution $f$ (\citeN{Peskir1}). Thus, to show uniqueness, it suffices to show that (\ref{eqn:unique}) is reducible to a Volterra equation of the second kind which has a unique solution.\\

Using Lemma (\ref{appendixlemma2}) we have, $|D_p(\frac{b(t)-b(s)}{\sqrt{t-s}})|<M_p$ for some $M_p>0$ and all $0\leq s\leq t,\ p\in \RR$.  Applying Abel's transform to equation (\ref{eqn:unique}) we obtain:
\begin{eqnarray}
\int_0^u\frac{g_{2\lambda+1}(t)}{(u-t)^{1-\lambda}}dt=\int_0^u\int_s^u\frac{K_{2\lambda+1}(t,s)}{(u-t)^{1-\lambda}(t-s)^{\lambda}}dtf(s)ds
\end{eqnarray}
where we have used Fubini's theorem (since $|K_{2\lambda+1}|$ is bounded) to exchange the order of integration. Let $\tilde{g}_{\lambda}(u)$ denote the left side of the above equation and $\tilde{K}_{2\lambda+1}(u,s)\triangleq\int_s^u\frac{K_{2\lambda+1}(t,s)}{(u-t)^{1-\lambda}(t-s)^{\lambda}}dt=\int_0^1\frac{K_{2\lambda+1}(y(u-s)+s,s)}{(1-y)^{1-\lambda}y^{\lambda}}dy$ then the last equation can be written as:
\begin{eqnarray}
\tilde{g}_{\lambda}(u)=\int_0^u\tilde{K}_{2\lambda+1}(u,s)f(s)ds \label{eqn:uniqueAbel}
\end{eqnarray}
Next, we apply the standard technique of differentiation on $u$ to reduce (\ref{eqn:uniqueAbel}) to a Volterra equation of the second kind. First we show that $g_{2\lambda+1}(0)=0$ and that $\tilde{g}_{\lambda}(u)$ has a continuous derivative for all $u\geq 0$. For any $\lambda_1,\lambda_2 \in \RR$, since $b(t)/\sqrt{t}\downarrow -\infty$ when $t\downarrow 0$ for a regular boundary $b$, we have (using the asymptotic expansion of the parabolic cylinder function)
\begin{eqnarray*}
\lim_{t\downarrow 0}\frac{g_{2\lambda_1+1}(t)}{t^{\lambda_2}}&=&\lim_{t\downarrow 0}\frac{e^{-b^2(t)/(2t)}}{t^{\lambda_1+\lambda_2}}\left(-\frac{b(t)}{\sqrt{t}}\right)^{2\lambda_1+1}\\
&=&\lim_{t\downarrow 0}e^{-\frac{b^2(t)}{2}\frac{1}{t}}(-b(t))^{2\lambda_1+1}\left(\frac{1}{t}\right)^{2\lambda_1+\lambda_2+1/2}=0
\end{eqnarray*}
since $\infty>-b(0)>0$. In particular $g_{2\lambda+1}(0)=0$. Also, using (\ref{eqn:pcfproperty3}),
\begin{eqnarray}
dg_{2\lambda+1}(t)/dt&=&-\frac{\lambda g_{2\lambda+1}(t)}{t}+\sqrt{t}g_{2\lambda+2}(t)\left(\frac{b'(t)}{\sqrt{t}}-\frac{b(t)}{2t^{3/2}}\right)\\
&=&-\frac{\lambda g_{2\lambda+1}(t)}{t}+\frac{g_{2\lambda+2}(t)}{t^{\epsilon}}b'(t)t^{\epsilon}-\frac{b(t)g_{2\lambda+2}(t)}{2t}
\end{eqnarray}
Under our assumption on the boundary and since $b(0)>-\infty$ each term in the last line goes to 0 as $t\downarrow 0$ and we obtain $$\lim_{t \downarrow 0}g'_{2\lambda+1}(t)=0.$$ Therefore, since $b$ is continuously differentiable, it follows that $g'_{2\lambda+1}(t)$ and $g_{2\lambda+1}(t)$ are continuous functions for all $t\geq 0$ and since $g_{2\lambda+1}(0)=0$, by Theorem 3, p.5, \citeN{Bocher09}, $\tilde{g}_{\lambda}(u)$ has a continuous derivative, for all $u\geq 0$, given by $$\tilde{g}'_{\lambda}(u)=\int_0^u\frac{g'_{2\lambda+1}(t)}{(u-t)^{1-\lambda}}dt$$\\

Next we compute the derivative, w.r.t. $u$, of the righthand side of (\ref{eqn:uniqueAbel}). Since $|D_{2\lambda+1}((b(s)-b(t))/\sqrt{t-s})|<C_{\lambda}$ for some $C_{\lambda}>0$ and all $0\leq s\leq t$ it follows that $K_{2\lambda+1}(y(u-s)+s,s)<C_{\lambda}$ for all $0\leq s\leq u$ and $0\leq y \leq 1$, while $\int_0^1C_{\lambda}\frac{1}{(1-y)^{1-\lambda}y^{\lambda}}dy=C_{\lambda}B(1-\lambda,\lambda)$ -- here $B(\cdot,\cdot)$ represents the Beta function. Thus, by the dominated convergence theorem,:
\begin{eqnarray*}
\tilde{K}_{2\lambda+1}(u,u)&=&\lim_{s\uparrow u}\tilde{K}_{2\lambda+1}(u,s)\\
&=&\int_0^1\lim_{s\uparrow u}K_{2\lambda+1}(y(u-s)+s,s)\frac{1}{(1-y)^{1-\lambda}y^{\lambda}}dy\\
&=&D_{2\lambda+1}(0)B(1-\lambda,\lambda)\neq 0
\end{eqnarray*}
for all $u\geq 0$.

Furthermore, using Lemma (\ref{appendixlemma2}),
\begin{align*}
y^{1/2-\epsilon}& |dK_{2\lambda+1}(y(u-s)+s,s)/du|\\
=&\left|\frac{e^{-\frac{(b(s)-b(y(u-s)+s))^2}{4y(u-s)}}D_{2\lambda+2}(\frac{b(s)-b(y(u-s)+s)}{\sqrt{y(u-s)}})}{(u-s)^{\epsilon+1/2}}(\frac{y(u-s)}{y(u-s)+s})^{\epsilon}\times
\right.\\
&\times \left. (b'(y(u-s)+s)(y(u-s)+s)^{\epsilon}-\frac{(b(y(u-s)+s)-b(s))(y(u-s)+s)^{\epsilon}}{2y(u-s)})\right|\\
:=&\frac{|H^{\epsilon}_{\lambda}(y,u,s)|}{(u-s)^{\epsilon+1/2}} \leq\frac{M}{(u-s)^{\epsilon+1/2}}
\end{align*}
for some constant $M>0$ and for all $0\leq s \leq u$ and $0\leq y\leq 1$. Also
\begin{align*}
&\int_0^u\left(\int_0^1\frac{M}{(u-s)^{\epsilon+1/2}}\frac{1}{(1-y)^{1-\lambda}y^{\lambda-(1/2-\epsilon)}}dy\right)f(s)ds\\
& \hspace{2cm}=MB(1-\lambda,\lambda-(1/2-\epsilon))\int_0^u\frac{f(s)}{(u-s)^{1/2+\epsilon}}ds<\infty
\end{align*}
since $f$ is continuous on $[0,T]$. Thus the derivative $dK_{2\lambda+1}(y(u-s)+s,s)/du$ is dominated by an integrable function. Denote $$K^{\delta}_{2\lambda+1}(y,u,s)\triangleq K_{2\lambda+1}(y(u-s)+s+\delta y,s)-K_{2\lambda+1}(y(u-s)+s,s)$$ for $\delta>0$ and note that $K^{\delta}_{2\lambda+1}(y,u,s)/\delta\rightarrow dK_{2\lambda+1}(y(u-s)+s,s)/du$ uniformly on $(s,y)\in [0,u]\times [0,1]$ as $\delta \downarrow 0$. Therefore, if $\mu$ denotes the measure on $[0,u]$ with Radon-Nykodim derivative $f$ and $\nu$ denotes the measure on $[0,1]$ with Radon-Nykodim derivative $1\left/(1-y)^{1-\lambda}y^{\lambda-(1/2-\epsilon)}\right.$, by Fubini's (applied twice) and the dominated convergence theorems, we have
\begin{align*}
\lim_{\delta \downarrow 0} & \int_0^u\int_0^1y^{1/2-\epsilon}\frac{K^{\delta}_{2\lambda+1}(y,u,s)}{\delta}\frac{f(s)}{(1-y)^{1-\lambda}y^{\lambda-(1/2-\epsilon)}}dyds\\
=&\lim_{\delta \downarrow 0}\int_{[0,u]\times[0,1]}y^{1/2-\epsilon}\frac{K^{\delta}_{2\lambda+1}(y,u,s)}{\delta}d(\mu \times \nu)\\
=&\int_{[0,u]\times[0,1]}\lim_{\delta \downarrow 0}y^{1/2-\epsilon}\frac{K^{\delta}_{2\lambda+1}(y,u,s)}{\delta}d(\mu \times \nu)\\
=&\int_{[0,u]\times[0,1]}\frac{H^{\epsilon}_{\lambda}(y,u,s)}{(u-s)^{\epsilon+1/2}}d(\mu \times \nu)\\
=&\int_0^u\left(\int_0^1\frac{H^{\epsilon}_{\lambda}(y,u,s)}{(1-y)^{1-\lambda}y^{\lambda-(1/2-\epsilon)}}dy\right)\frac{f(s)}{(u-s)^{\epsilon+1/2}}ds
\end{align*}
where the quantity $$ R^{\epsilon}_{\lambda}(u,s)\triangleq\int_0^1\frac{H^{\epsilon}_{\lambda}(y,u,s)}{(1-y)^{1-\lambda}y^{\lambda-(1/2-\epsilon)}}dy$$ is bounded by $MB(1-\lambda,\lambda-(1/2-\epsilon))$ and continuous for all $0\leq s<u$ with a possible discontinuity at $u=s$. Therefore the derivative of (\ref{eqn:uniqueAbel}) w.r.t. $u$ is given by
\begin{align*}
\tilde{g}'_{\lambda}(u)=&\lim_{\delta \downarrow 0}\left(\frac{\int_0^{u+\delta}\tilde{K}_{2\lambda+1}(u+\delta,s)f(s)ds-\int_0^u\tilde{K}_{2\lambda+1}(u,s)f(s)ds}{\delta}\right)\\
=&\lim_{\delta \downarrow 0}\frac{1}{\delta}\int_u^{u+\delta}\tilde{K}_{2\lambda+1}(u+\delta,s)f(s)ds\\
&+\lim_{\delta \downarrow 0}\int_0^u\frac{\tilde{K}_{2\lambda+1}(u+\delta,s)-\tilde{K}_{2\lambda+1}(u,s)}{\delta}f(s)ds\\
=&\tilde{K}_{2\lambda+1}(u,u)f(u)+ \int_0^u\frac{R^{\epsilon}_{\lambda}(u,s)}{(u-s)^{\epsilon+1/2}}f(s)ds
\end{align*}
Thus we obtained the Volterra equation of the second kind:
\begin{eqnarray}
\frac{\tilde{g}'_{\lambda}(u)}{\tilde{K}_{2\lambda+1}(u,u)}=f(u)+ \int_0^u\frac{R^{\epsilon}_{\lambda}(u,s)}{\tilde{K}_{2\lambda+1}(u,u)(u-s)^{\epsilon+1/2}}f(s)ds \label{voltera2unique}
\end{eqnarray}
Since $R^{\epsilon}_{\lambda}$ is finite on $0\leq s\leq u$ with a possible discontinuity only along the curve $s=u$ and since $f$ is continuous on $[0,\infty)$, by Theorem 3, p. 19, \citeN{Bocher09}, (\ref{voltera2unique}) has a unique continuous solution. Thus, the continuous solution to  (\ref{eqn:unique}) is unique. Therefore, for each $-1<p<1$ equation $B_p$ has a unique continuous solution. Furthermore, suppose $p\leq -1$ and $g:[0,t]\rightarrow \RR$ is any continuous solution to $B_p$. Then, using the integral representation of $D_p$ (and the fact that $D_p(x)>0,\ \forall p<0$) and Fubini's theorem, we can write equation $B_p$ as
\begin{eqnarray*}
\int_0^t\int_0^{\infty}e^{-\left(\frac{b(s)-b(t)}{\sqrt{t-s}}+\frac{u}{\sqrt{t-s}} \right)^2/2}u^{-p-1}\frac{g(s)}{\sqrt{t-s}}duds=\int_0^{\infty}\frac{u^{-p-1}}{\sqrt{t}}e^{-(u/\sqrt{t}-b(t))^2/(2t)}du\\
\int_0^{\infty}u^{-p-1}\int_0^{t}e^{-\left(\frac{b(s)-b(t)}{\sqrt{t-s}}+\frac{u}{\sqrt{t-s}} \right)^2/2}\frac{g(s)}{\sqrt{t-s}}dsdu=\int_0^{\infty}\frac{u^{-p-1}}{\sqrt{t}}e^{-(u/\sqrt{t}-b(t))^2/(2t)}du\\
\int_0^te^{-\left(\frac{b(s)-b(t)}{\sqrt{t-s}}+\frac{u}{\sqrt{t-s}} \right)^2/2}\frac{g(s)}{\sqrt{t-s}}ds=\frac{e^{-(u/\sqrt{t}-b(t))^2/(2t)}}{\sqrt{t}}
\end{eqnarray*}
where the last equality holds for all $u>0$ and follows from the uniqueness of the Mellin transform. Taking the limit $u\downarrow 0$ in the last equality it follows that any continuous solution to $\{B_p\}_{p\leq -1}$ is also a solution to $B_0$ which has a unique continuous solution. Thus we proved the following result:
\begin{theorem}
\label{thrm:unique}
For each $T>0$ let $b(t)$ be a regular boundary, continuously differentiable on $(0,T]$, and satisfy $|b'(t)|=O(t^{-\epsilon})$ for some $0<\epsilon<1/2$. Then $\tau$, the first-passage time of the standard Brownian motion $W_t$ to $b(t)$, has a continuous density function, $f$, given as the unique continuous solution of the class of integral equations $$\frac{e^{-\frac{b(t)^2}{4t}}D_{p}(-b(t)/\sqrt{t})}{t^{(p+1)/2}}=\int_{0}^te^{-\frac{(b(s)-b(t))^2}{4(t-s)}}\frac{D_{p}((b(s)-b(t))/\sqrt{t-s})}{(t-s)^{(p+1)/2}}f(s)ds $$ where $p<1$.
\end{theorem}

\subsection{Functional Transforms}
Next we consider some functional transforms of the boundary and the corresponding density functions. The new density functions can be easily expressed in terms of the original boundary and its density function using equation $B_0$ and Theorem \ref{thrm:unique}. Suppose $b$ satisfies the hypotheses of Theorem 3 with corresponding density function $f$ and introduce the functional transforms:
\begin{eqnarray}
(T_1^{\alpha}.b)(t)=b(t)+\alpha t, \ \alpha\in \RR\\
(T_2^{\gamma}.b)(t)=b(\gamma t)/\sqrt{\gamma}, \ \gamma>0\\
(T_3^{\beta}.b)(t)=(1+\beta t)b\left(\frac{t}{1+\beta t}\right),\ \beta \geq0
\end{eqnarray}
Note that we can set $\beta<0$ with $t\leq-1/\beta$ in the last transform. Moreover, $(T_1.b)$, $(T_2.b)$ and $(T_3.b)$ all satisfy the hypotheses of Theorem \ref{thrm:unique}. Denote with $f_1$, $f_2$ and $f_3$, respectively the corresponding density functions of the first-passage times of $W_t$ to these boundaries. Using equation $B_0$ we can easily find the relations between $f$ and $f_i,\ i=1,2,3$.

For $f_1$ we have:
\begin{align*}
\frac{e^{-(b(t)+\alpha t)^2/(2t)}}{\sqrt{t}}&=\int_0^t\frac{e^{-\frac{(b(t)-b(s)+\alpha (t-s))^2}{2(t-s)}}}{\sqrt{t-s}}f_1(s)ds\\
\Rightarrow \quad \frac{e^{-b(t)^2/(2t)}}{\sqrt{t}}&=\int_0^t\frac{e^{-\frac{(b(t)-b(s))^2}{2(t-s)}}}{\sqrt{t-s}}e^{\alpha b(s)+\alpha^2 s/2}f_1(s)ds
\end{align*}
Therefore, due to uniqueness of solutions (Theorem \ref{thrm:unique}), we must have
\begin{eqnarray}
f_1(t)=f(t)e^{-\alpha b(t)-\alpha^2 t/2}=e^{-\alpha(T^{\alpha/2}_1.b)(t)}f(t) \label {linbound}
\end{eqnarray}
This result can alternatively be obtained by a simple measure change argument.

For $f_2$ we obtain:
\begin{align*}
\frac{e^{-b(\gamma t)^2/(2\gamma t)}}{\sqrt{t}}&=\int_0^t\frac{e^{-\frac{(b(\gamma t)-b(\gamma s))^2}{2\gamma(t-s)}}}{\sqrt{t-s}}f_2(s)ds\\
\Rightarrow \quad \frac{e^{-b(u)^2/(2u)}}{\sqrt{u}}&=\int_0^u\frac{e^{-\frac{(b(u)-b(x))^2}{2(u-x)}}}{\sqrt{u-x}}f_2(x/\gamma)/\gamma dx
\end{align*}
and, applying Theorem \ref{thrm:unique},
\begin{eqnarray}
f_2(t)=\gamma f(\gamma t)=\gamma^{3/2}(T^{\gamma}_2.f)(t) \label{scalebound}
\end{eqnarray}
This result can alternatively be derived through the time change $t\rightarrow \alpha t$.

For $f_3$ we obtain:
\begin{align*}
&\frac{e^{-\frac{1}{2t}(1+\beta t)^2b\left(\frac{t}{1+\beta t}\right)^2}}{\sqrt{t}}=\int_0^t\exp\left(-\frac{1}{2(t-s)}\left((T^{\beta}_3.b))(t)-(T^{\beta}_3.b)(s)\right)^2\right)\frac{f_3(s)}{\sqrt{t-s}}ds\\
\Rightarrow &\frac{e^{-\frac{b^2(u)}{2u(1-\beta u)}}}{\sqrt{u}}=\int_0^u e^{-\frac{1-\beta x}{2(u-x)(1-\beta u)}\left(b(u)-b(x) \frac{1-\beta u}{1-\beta x}\right)^2} \frac{(1-\beta x)^{-3/2}f_3(\frac{x}{1-\beta x})}{\sqrt{u-x}}dx\\
\Rightarrow &\frac{e^{-\frac{b^2(u)}{2u(1-\beta u)}}}{\sqrt{u}}=e^{-\frac{\beta b^2(u)}{2(1-\beta u)}}\int_0^u e^{-\frac{(b(u)-b(x))^2}{2(u-x)}}e^{\frac{\beta b^2(x)}{1-\beta u}(\frac{1}{2}-\frac{\beta (u-x)}{2(1-\beta x)})} \frac{(1-\beta x)^{-3/2}f_3(\frac{x}{1-\beta x})}{\sqrt{u-x}}dx\\
\Rightarrow &\frac{e^{-b(u)^2/(2u)}}{\sqrt{u}}=\int_0^u\frac{e^{-\frac{(b(u)-b(x))^2}{2(u-x)}}}{\sqrt{u-x}}e^{\frac{\beta b^2(x)}{2(1-\beta x)}}(1-\beta x)^{-3/2}f_3\left(\frac{x}{1-\beta x}\right)dx
\end{align*}
where we have made the substitutions $u=\frac{t}{1+\beta t},\ x=\frac{s}{1+\beta s}$. Therefore, by Theorem \ref{thrm:unique},
\begin{eqnarray}
f_3(t)=f\left(\frac{t}{1+\beta t}\right)\exp\left(-\beta(1+\beta t) b^2(\frac{t}{1+\beta t})/2\right)(1+\beta t)^{-3/2} \label{patiealili}
\end{eqnarray}
This result was obtained by \citeN{AliliPatie2005} for more general boundaries, using probabilistic arguments.

It is instructive and pleasing that the integral equations lead to a unifying derivation of all of these transformation results. Combining the three transforms into the single transform $$(T.b)(t)\triangleq\frac{1+\beta \gamma t}{\sqrt{\gamma}}b\left(\frac{\gamma t}{1+\beta \gamma t}\right)+\alpha t,\ \alpha\in \RR,\ \beta \geq 0,\ \gamma >0$$ we have the following result for the corresponding density function $f_T$.
\begin{theorem}
For each $T>0$ let $b(t)$ be a regular boundary, continuously differentiable on $(0,T]$, and satisfy $|b'(t)|=O(t^{-\epsilon})$ for some $0<\epsilon<1/2$. Let $f$ be the continuous density function of the first-passage time of the standard Brownian motion $W_t$ to $b(t)$. Let $\tilde{b}(t)=(T.b)(t)$. Then the first-passage time to the boundary $\tilde{b}(t)$ has a continuous density $\tilde{f}$ given by:
\begin{eqnarray}
\tilde{f}(t)=\frac{\gamma f(\frac{\gamma t}{1+\beta \gamma t})}{(1+\beta \gamma t)^{3/2}}e^{-(1+\beta \gamma t)b(\frac{\gamma t}{1+\beta \gamma t})(\beta b(\frac{\gamma t}{1+\beta \gamma t})/2+\alpha/\sqrt{\gamma})-\alpha^2t/2} \label{funtransform}
\end{eqnarray}
\end{theorem}
\begin{proof}
 Since $b(t)$ is continuously differentiable then so is $\tilde{b}(t)$ and thus the first-passage time of $W_t$ to $\tilde{b}(t)$ has a continuous density function. Application of the transforms $T^{\beta}_3,\ T^{\gamma}_2,\ T^{\alpha}_1$ succesively to $b(t)$, transforms $f(t)$ to $\tilde{f}(t)$ and (\ref{funtransform}) follows from Theorem \ref{thrm:unique}. $\boxdot$
 \end{proof}

\section{Fredholm Equations}

Similarly to Section 2, in this section we examine the well known martingale $e^{-\alpha W_s-\alpha^2t/2}$ which gives rise to a Fredholm integral equation of the first kind. This equation is used to obtain alternative derivation of known closed form results for the linear, quadratic and square-root boundaries. Furthermore, as we will see, this equation is simply the Laplace transform of the integral equations (\ref{eqn:general}) for a particular class of boundaries. We assume that $b$ is continuous on $[0,\infty)$. Let $\tau_{\alpha}=\inf\left\{t>0;W_t \leq b(t)+\alpha t \right\}$ with cumulative distribution function $F_{\alpha}$ and define the set
\begin{eqnarray*}
\mathcal{A}_{b(t)}\triangleq\{ \alpha\in \RR;
b_{\alpha}(t)\triangleq b(t)+\alpha t \geq c \ for \ some \ c<0 \ and \ all \
t \geq 0 \}
\end{eqnarray*}
Under the measure $P^*$ given by $P^*(A)=\int_{A}Z(\omega)dP(\omega)$ where $Z=e^{-\frac{\alpha^2
t}{2}+\alpha W_t}$, $\tau_{\alpha}$ has distribution $F$. Then the
equality $E_{P^*}(1_{\tau_{\alpha} \leq t})=E_{P}(1_{\tau_{\alpha}
\leq t}Z)$ becomes
\begin{eqnarray*}
F(t)&=&\int_{0}^{t}E_{P}(e^{-\frac{\alpha^2 t}{2}+\alpha W_t}|\tau_{\alpha}=s)F_{\alpha}(ds)=\int_{0}^{t}E_{P}(e^{\alpha W_{\tau_{\alpha}}}e^{-\frac{\alpha^2 t}{2}+\alpha (W_t-W_s)}|\tau_{\alpha}=s)F_{\alpha}(ds)\\
&=&\int_{0}^{t}e^{\alpha (b(s)+\alpha s)}e^{-\frac{\alpha^2
t}{2}+\frac{\alpha^2 (t-s)}{2}}F_{\alpha}(ds) = \int_{0}^{t}e^{\alpha
b(s)+\frac{\alpha^2 s}{2}}F_{\alpha}(ds)
\end{eqnarray*}
where we have used the almost sure equality $W_{\tau_{\alpha}}=b(\tau_{\alpha})+\alpha \tau_{\alpha}$. Since the above is true for all $t \geq 0$ we have
\begin{eqnarray}
F_{\alpha}(dt)=F(dt)e^{-\alpha
b(t)-\frac{\alpha^2 t}{2}} \label{densities}
\end{eqnarray}
Under the assumption $b_{\alpha}(t) \geq c$ , we know
that $\tau_{\alpha} \leq \tau^c$ a.s. for all $\alpha \in \mathcal{A}$, where
$\tau^c\triangleq\inf\{t>0; W_t \leq c\}$. Since $\tau^c$ is almost surely
finite then so is $\tau_{\alpha}$, which implies $F_{\alpha}(\infty)=1$ and
thus, for $\alpha \in \mathcal{A}_{b(t)}$, using (\ref{densities}), we obtain the Fredholm
integral equation of the first kind
\begin{eqnarray}
\int_0^{\infty}e^{-\alpha b(s)-\frac{\alpha^2}{2}s}F(ds)=1 \label{eqn:mainreal}
\end{eqnarray}
with kernel $K(\alpha, s)=e^{-\alpha b(s)-\frac{\alpha^2s}{2}}$. Note that equation (\ref{densities}) holds for any boundary $b$ and $\alpha \in \RR$, while equation (\ref{eqn:mainreal}) holds for any $\alpha \in \mathcal{A}_{b(t)}$. The latter equation can also be derived using the martingale property of the Geometric Brownian motion together with the optional sampling theorem and has been found as early as \citeN{Shepp67}. \\
Next we extend equation (\ref{eqn:mainreal}) for complex values of $\alpha$. Consider the processes
\begin{eqnarray*}
X_t=e^{-xW_t+\frac{y^2-x^2}{2}t}\cos(yW_t+xyt)\\
Y_t=e^{-xW_t+\frac{y^2-x^2}{2}t}\sin(yW_t+xyt)
\end{eqnarray*}
for $x,y \in \RR$. Both processes are martingales for all real $x,y$
and thus the process $Z_t=X_t-iY_t=e^{-\alpha W_t-\frac{\alpha^2}{2}t}$
is a complex valued martingale where $\alpha=x+iy$. Define the class of continuous functions $b$.
$$\mathcal{B}=\{b(t); b(t)+ut>c,\ some\ c<0,\ \forall u<0,\ for\ large\ t \}$$ Note that
$b(t) \in \mathcal{B}$ implies $b(t)$ is uniformly bounded below and
thus the corresponding first passage time is almost surely finite.
\begin{theorem}
If $b \in \mathcal{B}$ and is continuous on $[0,\infty)$, then for all complex $\alpha$ with $|\arg(\alpha)|\leq \pi/2$,  the equality
\begin{eqnarray}
\int_0^{\infty}e^{-\alpha b(s)-\frac{\alpha^2}{2}s}F(ds)=1 \label{eqn:mainimaginary}
\end{eqnarray}
holds.
\end{theorem}
\begin{proof} First notice that for $b \in \mathcal{B}$ equation (\ref{eqn:mainreal}) holds
for all real $\alpha$ since $b(t)+\alpha t\geq b(t)-|\alpha|t>c$ for $t$ large enough. We first look at the quantity $E(e^{r\tau/2}),\ r>0$. For any such $r$, since $b\in \mathcal{B}$, there exists an $0<N(r)<\infty$ and a $c\in \mathbb R$ such that for $t>N(r)$ we have $\sqrt{r}b(t)>rt+c\sqrt{r}=\sqrt{r}(c+\sqrt{r}t)$. Then
\begin{eqnarray*}
E(e^{\frac{r}{2}\tau})&=&E(e^{\frac{r}{2}\tau}1(\tau\leq N(r)))+E(e^{\frac{r}{2}\tau}1(\tau>N(r)))\\
& \leq&E(e^{\frac{r}{2}\tau}1(\tau\leq N(r)))+E(e^{-c\sqrt{r}+\sqrt{r}b(\tau)-\frac{r}{2}\tau}1(\tau>N(r)))\\
&\leq&e^{\frac{r}{2}N(r)}+e^{-c\sqrt{r}}\int_{N(r)}^{\infty}e^{\sqrt{r}b(t)-rt/2}F(dt)\\
&\leq&e^{\frac{r}{2}N(r)}+e^{-c\sqrt{r}}\int_{0}^{\infty}e^{\sqrt{r}b(t)-rt/2}F(dt)<\infty
\end{eqnarray*}
Next we apply the optional sampling theorem by showing $E(|Z_{\tau}|)<\infty$ and $\lim_{t\rightarrow
\infty}E(Z_t1_{\tau>t})=0 $. For $x\geq 0$ and using the finitenes of $E(e^{r\tau/2}),\ r>0$, we have
\begin{align*}
E(|X_{\tau}|) \leq E(e^{-xb(\tau)+\frac{y^2-x^2}{2}\tau})\leq e^{-xc'}E(e^{\frac{y^2-x^2}{2}\tau}) \leq \infty\ .
\end{align*}
where $c'$ is the uniform lower bound of $b$. Similarly, for $x>0$, we can find an $M(x,y)$ such that for $t>M(x,y)$ we have the inequality $xb(t)-\frac{y^2-x^2}{2}t>0$. Then, for $t>M(x,y)$ we have $|X_t|1(\tau>t)\leq e^{-xb(t)+\frac{y^2-x^2}{2}t}1(\tau>t)<1(\tau>t)$ and thus $\lim_{t\uparrow \infty}E(|X_t|1_{\tau>t})\leq\lim_{t\uparrow \infty}P(\tau>t)=0$ since $\tau$  is almost surely finite. For $x=0$ $$E\left(|X_t|1(\tau>t)\right)\leq  E\left(e^{y^2t/2}1(\tau>t)\right)\leq E\left(e^{y^2/2\tau}1(\tau>t)\right)=\int_t^{\infty}e^{y^2s/2}F(ds)<\infty$$ and thus $\lim_{t \uparrow \infty}E\left(|X_t|1(\tau>t)\right)<\lim_{t \uparrow \infty}\int_t^{\infty}e^{y^2s/2}F(ds)=0$.
 Thus, for all $x\geq 0$, by the optional sampling theorem,  $X_t$ and $\tau$
satisfy $E(X_{\tau})=X_0=1$. The same arguments applied to the
process $Y_t$ yield $E(Y_{\tau})=Y_0=0$. Thus
$E(Z_{\tau})=E(X_{\tau})-iE(Y_{\tau})=1$ and the proof is
completed. $\boxdot$
\end{proof}

The above result gives an extension of the Fredholm equation (\ref{eqn:mainreal})
for boundaries belonging to the class $\mathcal{B}$. When
$-\frac{\pi}{4}\leq arg(\alpha)\leq \frac{\pi}{4}$ it is sufficient
that $b(t)$ is uniformly bounded below for equation (\ref{eqn:mainimaginary}) to hold
since $y^2-x^2<0$.\\

Fredholm equations of the first kind are notoriously difficult
to solve (even in the case when there is a unique solution). The two
general cases in which explicit results are available  are equations
with kernels of the form $K(\alpha t)$ or $K(\alpha - t)$. In the
first case we can obtain the Mellin transform of the solution and in
the second the Laplace transform. Next we examine boundaries which
give rise to such kernels. The following results are well known, however, here we demonstrate that they all follow from equations (\ref{eqn:mainreal}) and (\ref{eqn:mainimaginary}) and illustrate their importance.

\textbf{Example 1}.\\
$b(t)=-a+bt, \ a>0,b>0$. Thus $\mathcal{A}=\{\alpha \geq -b\}$ and equation
(\ref{eqn:mainreal}) becomes
\begin{eqnarray*}
\int_0^{\infty}e^{-(\frac{\alpha^2}{2}+\alpha b)t}f(t)dt=e^{-\alpha a}
\end{eqnarray*}
If $\tilde{f}$ is the Laplace transform of $f$ then the above equation
reads $\tilde{f}(u)=e^{-b+\sqrt{b^2+2u}}$ and this is the Laplace transform of the
well known Bachelier-Levy formula:
\begin{eqnarray*}
f(t)=\frac{a}{\sqrt{2 \pi}t^{3/2}}e^{-(a-bt)^2/(2t)}
\end{eqnarray*}

\textbf{Example 2}.\\
$b(t)=p\sqrt{t}- q, \ q\geq 0, \ p \neq 0$. Then
$\mathcal{A}=\{\alpha \geq 0\}$ and equation (\ref{eqn:mainreal}) becomes
\begin{eqnarray*}
\int_0^{\infty}e^{-\alpha p
\sqrt{t}-\frac{\alpha^2}{2}t}f(t)dt=e^{-\alpha q}
\end{eqnarray*}
Multiplying both sides of the above equation by $\alpha^{x-1}, \
x>0$ and integrating $\alpha$ on $[0,\infty)$ we obtain
\begin{align*}
& \int_0^{\infty} \alpha^{x-1} \int_0^{\infty} e^{-\alpha p
\sqrt{t}-\frac{\alpha^2}{2}t}f(t)dtd\alpha = \int_0^{\infty}\alpha^{x-1}e^{-\alpha
q}d\alpha\\
\Rightarrow \quad & \int_0^{\infty}f(t)dt\int_0^{\infty}\alpha^{x-1}e^{-\alpha p
\sqrt{t}-\frac{\alpha^2}{2}t}d\alpha = \frac{\Gamma(x)}{q^x}\\
\Rightarrow \quad & \int_0^{\infty}f(t)(2\frac{t}{2})^{-\frac{x}{2}}\Gamma(x)e^{p^2/4}D_{-x}(p)dt = \frac{\Gamma(x)}{q^x}\\
\Rightarrow \quad & \int_0^{\infty}t^{-\frac{x}{2}}f(t)dt = \frac{e^{-\frac{p^2}{4}}q^{-x}}{D_{-x}(p)}
\end{align*}
where $D$ is the parabolic cylinder function. The last equality gives us the Mellin transform of $f$ if we replace $x$ with $2(1-x), \ x<1$. Alternatively, by making the substitution $t=e^{u}$ in the last equation we obtain
\begin{eqnarray*}
\int_{-\infty}^{\infty}e^{-(\frac{x}{2}-1)u}f(e^{u})du=\frac{e^{-\frac{p^2}{4}}q^{-x}}{D_{-x}(p)}
\end{eqnarray*}
which gives us $\tilde{f}(x)$, the Laplace transform of $f(e^u)$, after
replacing $x$ with $2x+2$
\begin{eqnarray}
\tilde{f}(x)=\frac{e^{-\frac{p^2}{4}}q^{-(2x+2)}}{D_{-(2x+2)}(p)}, \ x>-1 \label{rootbound}
\end{eqnarray}
A similar approach was used in \citeN{Shepp67} for the first-passage time to the double boundary $\pm a\sqrt{t+b}$. \citeN{Novikov81} generalizes (\ref{rootbound}) for stable processes with a certain Laplace transform using a martingale approach.

\begin{figure}[t!]
\begin{center}
\includegraphics[width=6cm, height=7cm]{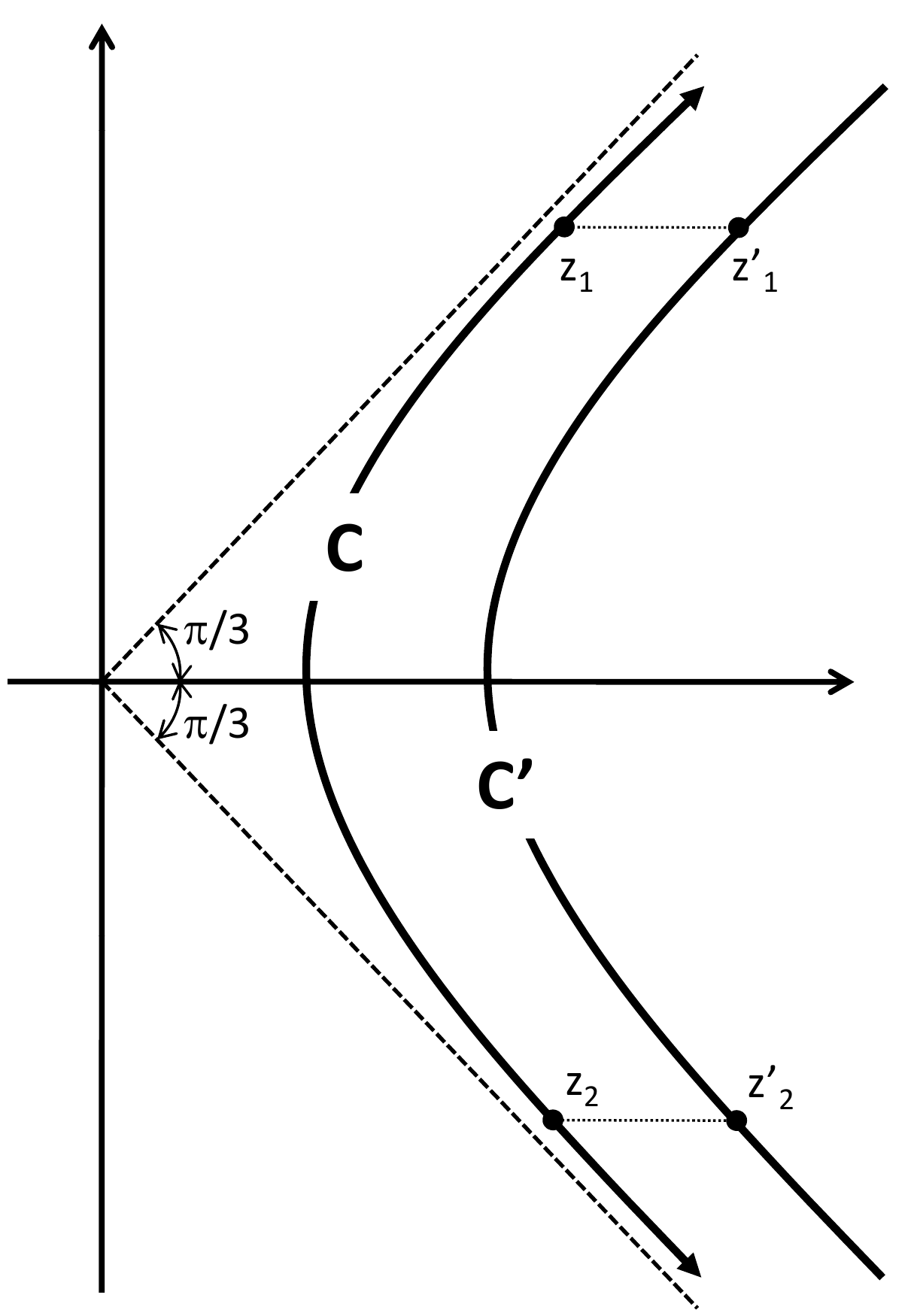}
\end{center}
\vspace{-1cm}
\caption{The contours of integration for quadratic boundaries in Example 3. \label{fig:contour}}
\end{figure}

\textbf{Example 3}.\\
$b(t)=\frac{pt^2}{2}-q, \ p,q>0$. Take $\alpha$ such that
$\Re(\alpha)>0$. Denote $\alpha'=\alpha(2p)^{1/3}$ then
$\Re(\alpha)>0$ as well.
 Using $\alpha'$ equation (\ref{eqn:mainimaginary}) becomes
\begin{eqnarray*}
\int_0^{\infty}e^{-\alpha'\frac{pt^2}{2}-\frac{\alpha'^2}{2}s}f(s)ds=e^{-\alpha'
q}
\end{eqnarray*}
and after completing the cube under the integral, multiplying both
sides of the equation by $\tfrac{e^{\alpha' \beta}}{2 \pi i}$,
$\beta>0$ and integrating $\alpha$ over any contour $C(\alpha)$ with end
points $\infty e^{-i\frac{\pi}{3}}$ and $\infty e^{i\frac{\pi}{3}}$
and $|\arg(\alpha)|\leq \pi/3$ (see Figure \ref{fig:contour}), we obtain
\begin{align*}
&\frac{1}{2\pi i}\int_C e^{\alpha' \beta}\int_0^{\infty}
e^{-\frac{(\alpha'+pt)^3}{6p}}e^{\frac{p^2t^3}{6}}f(t)dt d
\alpha = \frac{1}{2\pi i}\int_C e^{\alpha'(\beta-q)-\frac{\alpha'^3}{6p}}d\alpha\\
\Rightarrow & \int_0^{\infty}e^{-\beta pt}e^{\frac{p^2t^3}{6}}f(t)\int_C
e^{\beta(2p)^{1/3}(\alpha+\frac{pt}{(2p)^{1/3}})}e^{-\frac{1}{3}(\alpha+\frac{pt}{(2p)^{1/3}})^3}d
\alpha dt = \int_C
e^{\alpha(2p)^{1/3}(\beta-q)-\frac{\alpha^3}{3}}d \alpha
\end{align*}
The right hand side of the last equation is
$Ai((2p)^{1/3}(\beta-q))$, where $Ai$ is the Airy function (see (\ref{airy})). Next we
examine the contour integral on the left side. Define the contour
$C'=C+\frac{pt}{(2p)^{1/3}}$ and let $z_1$, $z_2$ be points on $C$
and their corresponding images on $C'$ be $z_1'$ and $z_2'$ (see Figure \ref{fig:contour}). Since
the function under the contour integral is analytic, its integral
over the simple closed contour $z_1z_1'z_2'z_2$ is $0$. Thus,
sending $z_1$ to $\infty e^{i\frac{\pi}{3}}$ and $z_2$ to $\infty
e^{-i\frac{\pi}{3}}$ we obtain
$$\int_C
e^{\beta(2p)^{1/3}(\alpha+\frac{pt}{(2p)^{1/3}})}e^{-\frac{1}{3}(\alpha+\frac{pt}{(2p)^{1/3}})^3}d
\alpha = \int_{C'}
e^{\beta(2p)^{1/3}\alpha}e^{-\frac{1}{3}\alpha^3}d \alpha =
Ai(\beta(2p)^{1/3})$$ since the contributions on the legs $z_1z_1'$ and
$z_2z_2'$ diminish in the limit. Therefore the Laplace transform of
$e^{\frac{p^2t^3}{6}}f(t)$ is given by
\begin{eqnarray}
\psi(\sigma)\triangleq\int_0^{\infty}e^{-\sigma
t}e^{\frac{p^2t^3}{6}}f(t)dt=\frac{Ai\left(\frac{2^{1/3}}{p^{2/3}}(\sigma-pq)\right)}{Ai(\sigma\frac{2^{1/3}}{p^{2/3}})} \label{quadraticbound}
\end{eqnarray}
The last example is an alternative derivation of (\ref{quadraticbound}) which was first obtained by \citeN{Salminen88} using measure change, and independently by \citeN{Groeneboom89} using a factorization of the density $f(t)$, involving a Bessel bridge and a killed Brownian motion.

\subsection{Connections with Volterra Integral Equations}

Finally, we discuss the connection between the Volterra integral equations of Section 2 (see \eqref{eqn:general})
%(with $p<1$)
and the Fredhom integral equations studied in the previous section (see \eqref{eqn:mainreal}) for a certain class of boundaries.

Let $b(t)$ be continuous and uniformly bounded below,i.e. there exists a constant $c<0$ such that for all $t\geq 0,\ b(t)>c$. Such boundaries satisfy (\ref{eqn:mainreal}) for all $\alpha \geq 0$. Set $y\leq c<0$ and $\alpha=\sqrt{2\beta},\ \beta \geq0$. Then, multiplying both sides of (\ref{eqn:mainreal}) by $\sqrt{\pi}2^{p+1/2}\beta^{p}e^{y\sqrt{2\beta}},\ p<0$ we obtain the equation:
\begin{eqnarray}
\int_0^{\infty}e^{-\beta s}\sqrt{\pi}2^{p+1/2}\beta^{p}e^{-\sqrt{2\beta}(b(s)-y)}F(ds)=\sqrt{\pi}2^{p+1/2}\beta^{p}e^{y\sqrt{2\beta}} \label{eqn:equivalence}
\end{eqnarray}
For any $-p,z,\beta>0,$ we have, (see \citeN{GradshteynRyzhik00}, 7.728),$$\int_0^{\infty}\frac{e^{-\beta x}x^{-(p+1)}}{\sqrt{\pi}2^{p+1/2}}e^{-\frac{z^2}{4x}}D_{2p+1}(z/\sqrt{x})dx=\beta^pe^{-z\sqrt{2\beta}}$$ and therefore the integral equation (\ref{eqn:equivalence}) can be written as:
\begin{align*}
&\int_0^{\infty}\int_s^{\infty}e^{-\beta t}\frac{e^{-\frac{(b(s)-y)^2}{4(t-s)}}D_{2p+1}(\frac{b(s)-y}{\sqrt{t-s}})}{(t-s)^{p+1}}dtF(ds)=\int_0^{\infty}e^{-\beta t}\frac{e^{-\frac{(y)^2}{4t}}D_{2p+1}(\frac{-y}{\sqrt{t}})}{t^{p+1}}dt\\
\Rightarrow &\int_0^{\infty}e^{-\beta t}\left[\int_0^t\frac{e^{-\frac{(b(s)-y)^2}{4(t-s)}}D_{2p+1}(\frac{b(s)-y}{\sqrt{t-s}})}{(t-s)^{p+1}}F(ds)\right]dt=\int_0^{\infty}e^{-\beta t}\left[\frac{e^{-\frac{(y)^2}{4t}}D_{2p+1}(\frac{-y}{\sqrt{t}})}{t^{p+1}}\right]dt
\end{align*}
Substituting $2p+1$ with $p$, and invoking the uniqueness of Laplace transforms, allows us to identify the terms in square braces and results in the class of integral equations
\begin{eqnarray}
\int_0^t\frac{e^{-\frac{(b(s)-y)^2}{4(t-s)}}D_{p}(\frac{b(s)-y}{\sqrt{t-s}})}{(t-s)^{(p+1)/2}}F(ds)=\frac{e^{-\frac{(y)^2}{4t}}D_{p}(\frac{-y}{\sqrt{t}})}{t^{(p+1)/2}},\ p<1 \label{eqn:mainrealtransform}
\end{eqnarray}
for all $0\leq s\leq t,\ p<1,\ y\leq c$. This is precisely the class of equations (\ref{eqn:general}) for $ p<1,\ y\leq c$, for boundaries $b(t)>c,\ t \geq0$.

\section{Conclusions}

In the first part of this article, we developed a new class of Volterra integral equations for the distribution of the first passage time(FPT)  of a standard Brownian motion to a regular boundary. This new class generalizes and unifies the class of all such previously known integral equations. Interestingly, this class arises through the optional stopping theorem applied to an interesting and new class of martingales generated by the parabolic cylinder functions. Through an Abel transformation, we were able to prove uniqueness of a solution to the integral equations. Based on uniqueness, we were then able to consolidate the derivation of the FPT distribution to a transformed boundary.

In the second part of this article, we generalized a class of Fredholm integral equations to the complex domain. These equations were then shown to provide a unified approach for computing the FPT distribution for linear, square root and quadratic boundaries. We believe that the method can be more widely applied by searching for specific factorizations of the kernel that produce known transforms such as Mellin, Laplace, Hilbert and so on. Finally we demonstrated that there is a fundamental connection between the Volterra and the Fredholm integral equations studied in this work.

There are several directions remaining open for future research.
\begin{itemize}
\item The first is clear but difficult: how can this larger (uncountably infinite) class of Volterra integral equations be used to extract the FPT distribution? One way is to explore the flexibility of the parabolic cylinder function and its connection to other special functions. Furthermore, the continuum of Volterra equations provides more flexibility for manipulation such as integration and differentiation w.r.t. the parameter $p$.

\item The search for new Volterra equations of the first kind is related to identifying analytical solution to the heat equation. The search for such solutions, which generate kernel functions with known properties, is another topic for future research. Any linear combination of solutions to the heat equation is also a solution and possibly, in the limit, one can obtain Volterra equations with more informative kernels.

\item We saw that taking the limit $y\uparrow b(t)$ in (\ref{eqn:general}), with $p=1$, produced the Volterra equation of the second kind (\ref{2kindvolterra}). This motivates the investigation of this limit for the equations with $p>1$. We suspect that in the computation of this limit for $p>1$, we can obtain new Volterra equations of the second kind and such equations are known to exhibit unique solutions and are generally easier to deal with than the Volterra equations of the first kind. However, such equations would hold for a restricted class of boundary functions.

\item The class of Volterra equations is also a useful tool for the inverse first passage time probem. Though, in this context, the equations are highly non-linear the generalization of this class provides flexibility for their manipulation which could extract new information.
\end{itemize}

\appendix

\section{Useful Bounding Lemmas}

In this appendix we provide two useful Lemmas which provide important bounds on the boundary and the density.
\begin{lemma}
\label{appendixlemma}
Suppose $b:(0,T]\rightarrow \RR$ is an increasing continuous function on $(0,\epsilon]$ for some $0<\epsilon<1$ with $b(0)=-\infty$. Let $h:\RR^+\rightarrow \RR$ and $h(x)=O(e^{ax^2})$ for large $x>0$ and some $0<a<1/2$. Define the first passage time $\tau\triangleq\{s>0;W_s\leq b(s)\}$. Then
\begin{eqnarray}
\int_0^{\epsilon}|h(-b(s))|F(ds)<\infty
\end{eqnarray}
where $F$ is the distribution of $\tau$.
\end{lemma}
\begin{proof} Without loss of generality we can assume $-b(t)\gg 0$ for $t\leq \epsilon$. Define the first-passage time $\tau_{b(s)}\triangleq\{t>0;W_t\leq b(s)\}$ for fixed $s<\epsilon$. Since $b(t)<b(s)$ for $t<s<\epsilon$ then $F(s)<F_{\tau_{b(s)}}(s)=2\Phi(-b(s)/\sqrt{s})$ for all $s\leq \epsilon$. Let $g(s)=h(-b(s)),\ s<\epsilon$ and fix $s_1<\epsilon$ and $\delta>0$ be such that $k_{a,\delta}(s_1)\triangleq e^{ab^2(s_1)}-\delta>0$. Define $s_n$ such that $e^{ab^2(s_n)}=k_{a,\delta}(s_1)+n\delta$. Since $e^{ab^2(s)}$ is monotone decreasing on $(0,\epsilon)$ with $e^{ab^2(0)}=\infty$ then $s_n\downarrow 0$ is a monotone decreasing sequence. Let
$$g_{\delta}(s)=\delta\sum_{n=2}^{\infty}1(s\leq s_{n-1})+(k_{a,\delta}(s_1)+\delta)1(s\leq \epsilon)$$ Then $0<|g(s)|\leq Me^{ab^2(s)}\leq Mg_{\delta}(s)$, $s\in (0,\epsilon]$, for some $M>0$, and by the dominated convergence theorem and the definition of $g$, there exists an $\epsilon>0$ such that
\begin{align*}
\int_0^{\epsilon}|g(s)|F(ds)\leq& M\int_0^{\epsilon}g_{\delta}F(ds)=M\delta\sum_{n=2}^{\infty}F(s_{n-1})+Me^{ab^2(s_1)}F(\epsilon)\\
\leq& 2M\delta\sum_{n=1}^{\infty}\Psi(b(s_n)/\sqrt{s_n})+C
\leq\frac{2\sqrt{2}M\delta}{\sqrt{\pi}}\sum_{n=1}^{\infty}\phi(b(s_n)/\sqrt{s_n})+C\\
\leq & \frac{2M\delta}{\pi}\sum_{n=1}^{\infty}e^{-b^2(s_n)/2}+C
=\frac{2M\delta}{\pi}\sum_{n=1}^{\infty}(k_{a,\delta}(s_1)+n\delta)^{-1/(2a)}+C\\
<&\infty \label{eqn:finite}
\end{align*}
where $C=Me^{ab^2(s_1)}F(\epsilon)$. The second inequality on the second line holds since $\Psi(x)\leq \phi(x),\ x>0$ while the last inequality follows from $a<1/2$. This completes the proof. $\boxdot$
\end{proof}
In particular, for $y\in \mathbb \RR,\ k>0$, we have
\begin{eqnarray}
\int_0^{\epsilon}(y-b(s))^kF(ds)<\infty \label{eqn:finite}
\end{eqnarray}

\begin{lemma}
\label{appendixlemma2}
Let $b(t)$ be a continuously differentiable function on $(0,T]$ with $-\infty<b(0)<0$ satisfying $\lim_{t \downarrow 0}|b'(t)|t^{\epsilon}<\infty$ for some $0<\epsilon<1/2$. Then, for all $0\leq s\leq t$,:
 $$\left|\frac{b(t)-b(s)}{\sqrt{t-s}}\right|<C\qquad and \quad
\left|\frac{b(t)-b(s)}{t-s}\right|t^{\epsilon}<K$$
for some positive constants $C$ and $K$.
\end{lemma}
\begin{proof} Since $b$ is continuous on $[0,T]$ the results hold for all $0\leq s< t\leq T$. Since $b$ is differentiable on $(0,T]$ the results hold on the curve $T\geq s=t>0$. We only need to check the case $s=t=0$. For $s=0$ and $t\downarrow 0$ we have that $$\lim_{t\downarrow 0}\left|\frac{b(t)-b(0)}{\sqrt{t}}\right|=\lim_{t\downarrow 0}\left|\frac{b(t)-b(0)}{t}\right|t^{\epsilon}t^{1/2-\epsilon}=\lim_{t \downarrow 0}|b'(t)|t^{\epsilon}t^{1/2-\epsilon}=0$$ and $$\lim_{t\downarrow 0}\left|\frac{b(t)-b(0)}{t}\right|t^{\epsilon}=\lim_{t \downarrow 0}|b'(t)|t^{\epsilon}<\infty \ .$$
This completes the proof. $\boxdot$
\end{proof}

\section{Some Special Functions}

In this Appendix we collect some important and useful results for the parabolic cylinder and Airy functions. For additional properties of the parabolic cylinder and Airy functions see \citeN{Erdelyi54} %pp. 116-126,
and \citeN{GradshteynRyzhik00}.

\subsection{The Parabolic Cylinder Function}
\label{pcf}
The parabolic cylinder functions are solutions to the differential equation
\begin{eqnarray}
\frac{d^2u}{dz^2}+\left(p+1/2-\frac{z^2}{4}\right)u=0 \ .\label{diffeqn1}
\end{eqnarray}
The complete set of solutions are $u(z)=D_p(z),\ D_p(-z),\ D_{-p-1}(iz),\ D_{-p-1}(-iz)$.

They admit the following integral representation for $p<0$:
\begin{eqnarray}
D_p(z)=\frac{e^{-z^2/4}}{\Gamma(-p)}\int_0^{\infty}e^{-xz-x^2/2}x^{-p-1}dx \ .\label{pcfintegral}
\end{eqnarray}

Furthermore, for special sets of parameters $p$, they reduce to other more well known special functions:
\begin{eqnarray}
D_n(z)&=& 2^{-n/2}e^{-z^2/4}H_n\left(\frac{z}{\sqrt{2}}\right) \ , \label{pcfHn} \\
D_p(z)&=& 2^{1/4+p/2}W_{1/4+p/2,-1/4}\left(\frac{z^2}{2} \right)z^{-1/2} \ , \label{pcfW} \\
D_{-1/2}(z)&=& \sqrt{z\pi/2}K_{1/4}\left(\frac{z^2}{4}\right) \ , \label{pcfK} \\
D_{-2}(z)&=& e^{z^2/4}\left(e^{-z^2/2}-\sqrt{2\pi}z\Phi(-z) \right) \ . \label{Dm2}
\end{eqnarray}
Here, $H_n$ is the Hermite polynomial of degree $n$, $W$ is the Whittaker function, $K$ is the modified Bessel function of the third kind.

For large argument $z$, they admit the following asymptotic expansions:
\begin{eqnarray}
D_p(z) &\sim& \frac{\sqrt{2\pi}}{\Gamma(-p)}e^{i\pi p}z^{-p-1}e^{z^2/4},\ |z|\uparrow \infty,\ \pi/4<|\arg(z)|<5\pi/4  \ , \label{pcfassymptotic-}\\
D_p(z)&\sim& z^pe^{-z^2/4},\ |z|\uparrow \infty,\ |\arg(z)|<3\pi/4 \label{pcfassymptotic+} \ .
\end{eqnarray}

They are closed under derivative and integral operations in the following sense:
\begin{eqnarray}
\frac{d}{dz}D_p(z)&=&\frac{1}{2}zD_p(z)-D_{p+1}(z) \ ,\label{eqn:pcfproperty1}\\
\frac{d}{dz}e^{-z^2/4}D_p(z)&=&e^{-z^2/4}D_{p+1}(z) \ ,\label{eqn:pcfproperty3}\\
\int_0^{\infty}e^{-z^2/4}D_{-p}(z)dz&=&\frac{\sqrt{\pi} 2^{-p/2-1/2}}{\Gamma(p/2+1)}=D_{-(p+1)}(0) \ . \label{eqn:pcfproperty2}
\end{eqnarray}

\subsection{The Airy function}
The Airy function (for complext argument $x$) has the integral representation:
\begin{eqnarray}
Ai(x)=\frac{1}{2\pi i}\int_C e^{t^3/3-xt}dt \label{airy}
\end{eqnarray}
where the integral is over a path $C$ with end points $\infty e^{-\pi i/3}$ and $\infty e^{\pi i/3}$ (see Figure \ref{fig:contour}).

\bibliographystyle{chicago}
	\bibliography{FPT}

\end{document}